\documentclass[%
aip,
amsmath,amssymb,
reprint,%
]{revtex4-1}

\usepackage{graphicx}% Include figure files
\usepackage{dcolumn}% Align table columns on decimal point
\usepackage{bm}% bold math
\usepackage[utf8]{inputenc}
\usepackage[T1]{fontenc}
\usepackage{mathptmx}
\usepackage{etoolbox}

\usepackage{float}

\newtheorem{theorem}{Theorem}

\newtheorem{proposition}[theorem]{Proposition}

\newtheorem{example}[theorem]{Example}
\newtheorem{definition}[theorem]{Definition}
\newtheorem{problem}[theorem]{Problem}
\newtheorem{algorithm}[theorem]{Algorithm}

\newtheorem{remark}[theorem]{Remark}

\makeatletter
\def\@email#1#2{%
	\endgroup
	\patchcmd{\titleblock@produce}
	{\frontmatter@RRAPformat}
	{\frontmatter@RRAPformat{\produce@RRAP{*#1\href{mailto:#2}{#2}}}\frontmatter@RRAPformat}
	{}{}
}%
\makeatother

\usepackage{hyperref}
\usepackage{mathtools}
\mathtoolsset{showonlyrefs}
\begin{document}
	
	\preprint{AIP/123-QED}
	
	\title[Neural Network Approximation of Optimal Controls for Stochastic Reaction-Diffusion Equations]{Neural Network Approximation of Optimal Controls for Stochastic Reaction-Diffusion Equations}
	\author{W. Stannat}
	\author{A. Vogler}
    \affiliation{Institute of Mathematics, Technische Universität Berlin, Straße des 17. Juni 136, 
        10623 Berlin, Germany}
	\author{L. Wessels}
	\affiliation{School of Mathematics, Georgia Institute of Technology, 686 Cherry Street, Atlanta, 
	   GA 30332-0160, USA}
    \email{stannat@math.tu-berlin.de, vogler@math.tu-berlin.de, wessels@gatech.edu.}
	
	\date{14 September 2023}% It is always \today, today,
	
	\begin{abstract}
		We present a numerical algorithm that allows the approximation of optimal controls for stochastic reaction-diffusion equations with additive noise by first reducing the problem to controls of feedback form and then approximating the feedback function using finitely based approximations. Using structural assumptions on the finitely based approximations, rates for the approximation error of the cost can be obtained. Our algorithm significantly reduces the computational complexity of finding controls with asymptotically optimal cost. Numerical experiments using artificial neural networks as well as radial basis function networks illustrate the performance of our algorithm. Our approach can also be applied to stochastic control problems for high dimensional stochastic differential equations and more general stochastic partial differential equations.
	\end{abstract}
	
	\maketitle
	
	\begin{quotation}
		There is a huge body of literature on optimal control problems with partial differential equation (PDE) constraints and their numerical treatment. Recent years have shown a rising interest in the optimal control of stochastic partial differential equations (SPDEs). However, the numerical approximation of optimal controls, let alone its practical implementation in the stochastic case, faces serious obstacles due to the computational complexity of classical algorithms. In this work, we present a new numerical algorithm that approximates feedback controls for SPDEs with asymptotically optimal cost. The algorithm is based on adjoint calculus applied to gradient descent. For the approximation of feedback controls, we use finitely based approximations such as artificial neural networks or radial basis function networks. The restriction to additive noise and our approach for the approximation of feedback functions enables us to significantly reduce the algorithmic complexity of our approach in comparison with classical algorithms.
		%Please establish the first paragraph of the article as a Lead Paragraph containing the main points of the article in terms accessible to non-specialist readers, and provide the “big picture”.}
\end{quotation}

\section{Introduction}
\label{Introduction}

For a fixed finite time horizon $T>0$, we consider the randomly forced reaction-diffusion equation
\begin{equation}\label{pde}
	\begin{cases}
		\partial_t u(t,x) = \Delta u(t,x) + f(u(t,x)) + \xi(t,x),\, (t,x)\in [0,T]\times \Lambda\\
		u(0,x)= u(x),\quad x\in \Lambda\\
        \frac{\partial u}{\partial \textbf{n}}(t,x) = 0,\quad (t,x)\in (0,T]\times \partial \Lambda
	\end{cases}
\end{equation}
on a bounded domain $\Lambda\subset\mathbb{R}^d$. Here, $\Delta := \sum_{k=1}^d \frac{\partial^2}{\partial x^2_k}$ denotes the Laplace operator, $f:\mathbb{R}\to\mathbb{R}$ models a local reaction term, $\xi(t,x)$ are random fluctuations, and $\frac{\partial u}{\partial \textbf{n}}$ denotes the partial derivative of $u$ with respect to the unit outward normal of $\partial \Lambda$, $\textbf{n}$. 
Deterministic reaction-diffusion equations are ubiquitous in the natural sciences and, in many situations, taking into account random fluctuations leads to the more realistic model \eqref{pde}.
Typically, these random fluctuations are highly irregular, and therefore, in order to treat this equation rigorously, we reformulate it as the following $L^2(\Lambda)$-valued SPDE:
\begin{equation}\label{spde1}
	\begin{cases}
		\mathrm{d}u_t=[ \Delta u_t+\mathcal{F}(u_t) ]\mathrm{d}t+ \sigma \mathrm{d}W_t, \quad t\in [0,T]\\
		u_0=u\in L^2(\Lambda),
	\end{cases}
\end{equation}
where $u_t(x)=u(t,x)$ for fixed $t\in [0,T]$ is considered as an element in $L^2(\Lambda)$. The random fluctuations $\xi(t,x)$ are modeled by a cylindrical Wiener process $(W_t)_{t\in[0,T]}$ on $L^2(\Lambda)$, defined on some underlying probability space $(\Omega, \mathcal{F},\mathbb{P})$, and $\sigma:L^2(\Lambda)\to L^2(\Lambda)$ is a Hilbert-Schmidt operator. Note that the boundary conditions are encoded in the choice of the Gelfand triple $H^1(\Lambda) \hookrightarrow L^2(\Lambda) \hookrightarrow (H^1(\Lambda))^{\ast}$, where $H^1(\Lambda)$ denotes the Sobolev space of order $1$, see Example 4.D in Chapter 1 in \cite{showalter1997}. For more details on the mathematical theory of SPDEs, see \cite{liu2015}. Furthermore, $\mathcal{F}:L^2(\Lambda)\to L^2(\Lambda)$ denotes the Nemytskii operator associated with $f$, i.e.,
\begin{equation}
	\mathcal{F}(u)(x) := f(u(x)),\quad u\in L^2(\Lambda),\, x\in \Lambda.
\end{equation}

The objective of control theory is to achieve a desired outcome for a dynamical system by applying an external input which can be chosen freely among a set of admissible inputs. In order to set up a mathematical formulation for the control of the SPDE \eqref{spde1}, we introduce a control process $\mathfrak{g}:[0,T]\times \Omega \to \mathcal{U}\subset L^2(\Lambda)$, adapted to the filtration generated by $u_t$, as a forcing term on the right-hand side of the equation
\begin{equation}\label{intro:controlledstate}
	\begin{cases}
		\mathrm{d}u^{\mathfrak{g}}_t=[ \Delta u^{\mathfrak{g}}_t+\mathcal{F}(u^{\mathfrak{g}}_t)+\mathfrak{g}_t ]\mathrm{d}t+ \sigma \mathrm{d}W_t, \quad t\in [0,T]\\
		u^{\mathfrak{g}}_0=u\in L^2(\Lambda),
	\end{cases}
\end{equation}
and define the cost functional
\begin{multline}\label{intro:costfunctional}
	J(\mathfrak{g}):=\mathbb{E} \bigg [\int_{0}^{T}\int_{\Lambda} l(t,x,u^{\mathfrak{g}}_t(x)) \mathrm{d}x + \frac{\nu}{2} \| \mathfrak{g}_t \|^2_{L^2(\Lambda)} \mathrm{d}t\\
	+ \int_{\Lambda}m(x,u^{\mathfrak{g}}_T(x))\mathrm{d}x \bigg ].
\end{multline}
Here, the running cost $l (t,x,u) :[0,T]\times \Lambda\times \mathbb{R} \rightarrow \mathbb{R}$ and the terminal cost $m (x,u) : 
\Lambda\times\mathbb{R} \rightarrow \mathbb{R}$ are assumed to be differentiable and locally Lipschitz continuous in the state variable $u$. Now, the objective is to find a control process $\mathfrak{g}$ that minimizes the cost functional \eqref{intro:costfunctional} over some set of admissible controls subject to \eqref{intro:controlledstate}.

In the deterministic case, there is a huge body of literature on necessary and sufficient optimality conditions\cite{hinze2009,li1995,troeltzsch2010}, as well as numerical algorithms which efficiently approximate optimal controls, see e.g. \cite{buchholz2013,buchholz20132,ryll2016,ryll2017}. In recent years, extensions to the stochastic case have seen a rising interest in the mathematical literature leading to necessary and sufficient optimality conditions in great generality \cite{cordoni2018,du2013,fabbri2017,frankowska2020,fuhrman2013,fuhrman2018,fuhrman2016,lue2014,lue2015,lue2018,stannat2021,stannat2022,wessels2022}. However, classical algorithms for the numerical approximation of optimal controls in the stochastic case either require the approximation of backward SPDEs or the approximation of infinite dimensional Hamilton-Jacobi-Bellman (HJB) equations. Due to the curse of dimensionality, both of these alternatives are computationally very expensive, leading to an increased interest in the development of new, more efficient algorithms \cite{beck2021,beck2019,dolgov2021,dunst2019,dunst2016,e2017,gorodetsky2018,kalise2018,nuesken2021,oster2022,richter2021,sirignano2018}.

A considerable reduction of the problem can be achieved in the case when 
an optimal control $\mathfrak{g}^\ast$ is of feedback-type, i.e.,
\begin{equation}\label{feedbackIntro}
	\mathfrak{g}^{\ast}_t = G^{\ast} (t,u^{G^{\ast}}_t ),
\end{equation}
for some $G^{\ast}:[0,T]\times L^2(\Lambda) \to L^2(\Lambda)$. 
In this work, we consider a mathematical setting in which the optimal 
control is indeed of the above type. In a first step we then use finitely 
based approximations for $G^{\ast}$ in order to approximate the optimal 
control $\mathfrak{g}^{\ast}$. This approach together with the restriction 
to additive noise enables us to replace the backward SPDE arising in adjoint 
calculus by a random backward PDE which significantly reduces its 
computational complexity. A similar idea was already used in 
\cite{stannat20212} in the case of deterministic controls which enter the state equation linearly. However, in the case of feedback controls, the control does not enter the equation linearly, thus requiring an extension of the results in \cite{stannat20212}.

The reduction to finitely based feedback controls then allows in a second step the local uniform approximation of finitely based feedback controls in appropriate ansatz spaces of variable, but finite, dimension. Solving the corresponding finite dimensional optimal control problem, e.g. with artificial neural networks, leads to an efficient computation of finite dimensional controls, whose costs approximate the optimal cost with increasing dimension, see Theorem \ref{main1}. Imposing additional smoothness assumptions on the optimal control $\mathfrak{g}^{\ast}_t  = G^{\ast} (t,u^{G^{\ast}}_t )$, in particular $G^\ast$ globally Lipschitz, also enables us to derive rates for the approximation error in ansatz spaces consisting of Lipschitz continuous feedback functions, see Theorem \ref{main2}. We discuss the example of neural network approximation in detail; see Theorem \ref{thm:main3} for the universal approximation of finitely based feedback controls with one-layer artificial neural networks and Proposition \ref{neuralapprox} for rates on the approximation error.

The practical implementation of the gradient descent algorithm for the approximation of finite dimensional feedback controls also requires the numerical discretization of the controlled state equation. In Proposition \ref{prop:numericaldiscretization}, we derive the combined approximation error for the optimal cost. Detailed proofs of our results, that also hold for more general ansatz spaces beyond artificial neural networks, can be found in \cite{stannat2023}. The performance of our gradient descent algorithm is illustrated with three examples. The first example deals with the validation of our algorithm in the case of linear quadratic control. The other two examples consider the problem of stabilizing a bump solution of the stochastic Nagumo equation, in one case with general feedback controls, in the other case with feedback controls of Nemytskii-type. 

The remainder of the paper is organized as follows: First, in Section \ref{reduction}, we show that it is sufficient to consider feedback controls. In Section \ref{Ansatz}, we explain how to approximate the optimal control by introducing ansatz spaces of controls that are suited for the numerical implementation and present our main results. In Section \ref{gradientdescent}, we describe our gradient descent algorithm. Next, in Section \ref{examples}, we discuss the explicit example of artificial neural networks for the finite dimensional ansatz spaces. In Section \ref{numericaldiscretization}, we investigate the error resulting from the numerical discretization. Finally, in Section \ref{simulations}, we present numerical experiments using artificial neural networks and radial basis function networks.

\section{\label{reduction}Optimal Controls of Feedback-Type}

In order to reduce the complexity of the problem, we assume that there exists an optimal control $\mathfrak{g}^{\ast}$ in feedback form \eqref{feedbackIntro}, for some continuous $G^{\ast}:[0,T]\times L^2(\Lambda) \to L^2(\Lambda)$ that satisfies a linear growth condition. This can be achieved, using the solution $V:[0,T]\times L^2(\Lambda) \to \mathbb{R}$ of the associated HJB equation
\begin{equation}\label{HJB}
	\begin{cases}
		\partial_t V+\frac{1}{2}\text{tr}(\sigma^{\ast}D^2V \sigma)+\langle DV,\Delta u+\mathcal{F}(u)\rangle_{L^2(\Lambda)}\\
		+\int_{\Lambda}l(t,x,u(x))\mathrm{d}x
		+ \inf_{G\in \mathcal{U}}\{\langle DV,G\rangle+ \frac{\nu}{2} \|G\|_{L^2(\Lambda)}^2\} =0,\\
		\qquad\qquad\qquad\qquad\qquad\qquad\qquad (t,u) \in [0,T] \times L^2(\Lambda)\\
		V(T,u)=\int_{\Lambda}m(x,u(x))\mathrm{d}x,\quad u\in L^2(\Lambda),
	\end{cases}
\end{equation}
where $DV$ and $D^2V$ denote the first and second Fr\'echet derivative of $V$ with respect to the second variable. For solution theories regarding equations of this type, see \cite{fabbri2017}. In particular, if \eqref{HJB} has a unique mild solution satisfying certain regularity assumptions, and
\begin{equation}
	\gamma(p):=\text{arginf}_{G\in \mathcal{U}}\left \{\langle p,G\rangle+ \frac{\nu}{2} \|G\|_{L^2(\Lambda)}^2\right \}
\end{equation}
is continuous, then 
\begin{equation}
	\mathfrak{g}^{\ast}_t :=G^{\ast}(t,u_t^{G^{\ast}}):=\gamma(DV(t,u_t^{G^{\ast}})).
\end{equation}
is an optimal control in feedback form, provided that $u^{G^{\ast}}$ is a unique strong solution of the closed loop equation 
\begin{equation}\label{state}
	\begin{cases}
		\mathrm{d}u^G_t=[ \Delta u^G_t+\mathcal{F}(u^G_t)+G(t,u^G_t) ]\mathrm{d}t+ \sigma \mathrm{d}W_t, \; t\in [0,T]\\
		u^G_0=u\in L^2(\Lambda)
	\end{cases}
\end{equation}
with $G=G^{\ast}$. A sufficient condition for \eqref{state} to have a unique strong solution is that $G^{\ast}$ is Lipschitz continuous in $u$, see e.g. \cite{liu2015}. This is in particular the case if the solution $V$ of the HJB equation \eqref{HJB} has a bounded second derivative in $u$, see Theorem 4.155 and Remark 4.202 in \cite{fabbri2017}. However, directly tackling the optimal control problem by approximating the solution of the HJB equation \eqref{HJB} numerically is very challenging due to the infinite dimensionality of the domain of $V$. Instead, our approach is to approximate the feedback function directly. Therefore we consider the following feedback control problem: Minimize
\begin{multline}\label{cost}
	J(G):=\mathbb{E} \bigg [\int_{0}^{T}\int_{\Lambda} l(t,x,u^G_t(x)) \mathrm{d}x + \frac{\nu}{2} \| G(t,u^G_t) \|^2_{L^2(\Lambda)} \mathrm{d}t\\
	+ \int_{\Lambda}m(x,u^G_T(x))\mathrm{d}x \bigg]
\end{multline}
subject to equation \eqref{state}. We seek to minimize $J$ over the set of admissible controls
\begin{equation}
	U_{\text{ad}} := \{ G :[0,T]\times L^2(\Lambda)\to \mathcal{U} |\, G(\cdot,u_{\cdot}^G)\in \mathbb{A}  \},
\end{equation}
where $\mathbb{A} = L^2([0,T]\times\Omega;L^2(\Lambda))$.

\begin{example}[Linear Quadratic Control]\label{linearquadratic}
	
	Let us consider the linear quadratic control problem
	\begin{equation}
		\begin{cases}
			\mathrm{d}u^{\mathfrak{g}}_t=[ \Delta u^{\mathfrak{g}}_t + \mathfrak{g}_t ]\mathrm{d}t+ \sigma \mathrm{d}W_t, \quad t\in [0,T]\\
			u^{\mathfrak{g}}_0=u\in L^2(\Lambda),
		\end{cases}
	\end{equation}
	where $\mathfrak{g}:[0,T]\times \Omega \to L^2(\Lambda)$ is an adapted process, and
	\begin{equation}
		J(\mathfrak{g}):=\frac12 \mathbb{E}\left[\int_{0}^{T}\int_{\Lambda} (u^{\mathfrak{g}}_t(x))^2 + \mathfrak{g}_t^2(x) \mathrm{d}x \mathrm{d}t + \int_{\Lambda} (u^{\mathfrak{g}}_T(x))^2\mathrm{d}x \right].
	\end{equation}
	In this case, the optimal control $\mathfrak{g}^{\ast}$ is indeed of feedback form, given by
	\begin{equation}
		\mathfrak{g}^{\ast}_t = P(t) u^{\mathfrak{g}^{\ast}}_t
	\end{equation}
	where $P:[0,T]\to L(L^2(\Lambda))$ is the solution of the Riccati equation \cite{tudor1990}
	\begin{equation}\label{Ricc}
		\begin{cases}
			\partial_t P(t) + P(t) \Delta +\Delta P(t) - \text{Id} + P^2(t) = 0,\;\; t\in [0,T]\\
			P(T) = - \text{Id}.
		\end{cases}
	\end{equation}
\end{example}

\section{Construction of Ansatz Spaces and Main Results}\label{Ansatz}

For the efficient numerical implementation we need to restrict ourselves to a subset of controls $\mathbb{U}\subset U_{\text{ad}}$ that are suitable for the implementation in a gradient descent algorithm. However, we need to ensure that when minimizing over $\mathbb{U}$, we do not end up much worse than in the original control problem. In this section we will provide a method to construct suitable ansatz spaces $\mathbb{U}$ when we do not have any particular control constraints, i.e. $\mathcal{U}=L^2(\Lambda)$.

In order to construct a suitable space $\mathbb{U}$, we introduce the so-called finitely based approximation of a function $G:[0,T]\times L^2(\Lambda)\rightarrow L^2(\Lambda)$. To this end, we consider finite dimensional subspaces 
\begin{align}
	S_n\subset L^2(\Lambda)
\end{align}
with orthonormal basis $e_0,\dots,e_n$ and orthogonal projections $P_n:L^2(\Lambda)\rightarrow S_n$, such that
\begin{align}
	\|P_nu-u\|_{L^2(\Lambda)}\rightarrow 0.
\end{align}
The finitely based approximations of $G$ with respect to $S_n$ are then defined by
\begin{align}\label{finitelybased}
	G^n(t,u) &:= P_n G(t,P_n u)\\
	&=\sum_{k=0}^n\langle G(t,\sum_{j=0}^n u_j e_j),e_k\rangle_{L^2(\Lambda)}e_k,
\end{align}
where $u_j:=\langle u,e_j\rangle_{L^2(\Lambda)}$. One possible choice for the finite dimensional subspaces in the case of $\Lambda = (0,1)$ is
\begin{equation}
	S_n:=\text{span}\left \{ 1,\sqrt{2}\cos(k \pi\,\cdot\,)\middle |\,k=1,\dots,n \right \}.
\end{equation}
In this case the finitely based approximation of $G$ is given by
\begin{multline}
	G^n(t,u)(x) \\:= \sum_{k=0}^{n}2\langle G(t,\sum_{j=0}^{n}u_j\sqrt{2}\cos(j\pi\,\cdot\,)),\cos(k\pi\,\cdot\,)\rangle_{L^2(0,1)}\cos(k\pi x),
\end{multline}
for $u=\sum_{j=0}^\infty u_j \cos(j\pi\,\cdot\,)\in L^2(\Lambda)$.

Our main results provide approximation results for ansatz spaces $\mathbb{U}\subset U_{\text{ad}}$. The main assumption on the ansatz spaces $\mathbb{U}\subset U_{\text{ad}}$ is the following approximation property with respect to the optimal feedback $G^{\ast}$:

\begin{definition}\label{def:universal}
	Let $G:[0,T]\times L^2(\Lambda)\rightarrow L^2(\Lambda)$. We say that a subset $\mathbb{U}\subset U_{\text{ad}}$ satisfies the uniform approximation property with respect to $G$, if there exists a sequence $(G^{n,m})_{n,m\in \mathbb{N}}\subset \mathbb{U}$ that satisfies for any $n\in \mathbb{N}$ the following linear growth condition uniformly in $m$
	\begin{align}
		\|G^{n,m}(t,u)\|_{L^2(\Lambda)}&\leq C_n(1+\|u\|_{L^2(\Lambda)}),
	\end{align}
	for some constant $C_n>0$, such that for any $R>0$ and $n\in \mathbb{N}$
	\begin{align}
		\lim\limits_{m\rightarrow \infty}\sup_{(t,u)\in [0,T]\times \mathcal{B}_{L^2(\Lambda)}(0,R)}\|G^n(t,u)-G^{n,m}(t,u)\|_{L^2(\Lambda)}^2=0,
	\end{align}
	where 
	\begin{align}
		\mathcal{B}_{L^2(\Lambda)}(0,R):=\{u\in L^2(\Lambda)|\|u\|_{L^2(\Lambda)}&\leq R\}
	\end{align}
	and $G^n$ is given as in \eqref{finitelybased}.
\end{definition}
We would like to stress that in the above definition we only require that the finitely based approximations of $G$ can be approximated uniformly on bounded sets, and not $G$ itself.

In the first part of this section we will explain how to construct ansatz spaces $\mathbb{U}$ that satisfy the uniform approximation property with respect to the optimal feedback $G^{\ast}$ and present our main approximation result for this type of ansatz space. These ansatz spaces consider finitely based controls of arbitrary dimension $ n\in \mathbb{N}$, which are of the type 
\begin{align}
	G(t,u)=\sum_{k=0}^n \psi^k(t,P_nu)e_k,
\end{align}
for some functions $\psi^k:[0,T]\times S_n\rightarrow \mathbb{R}$. However, in practice the dimension $n\in \mathbb{N}$ of the ansatz space needs to be fixed a priori and therefore one is interested in how close one can get to the optimal cost. Our second main result, Theorem \ref{main2}, provides explicit convergence rates, however, we need to strengthen the assumption on our ansatz spaces. In the second part of this section we will explain how to construct for any $n\in \mathbb{N}$ a sequence of ansatz spaces $(\mathbb{U}^{n,m})_{m\in \mathbb{N}}\subset U_{\text{ad}}$ that satisfies the following approximation property with respect to the optimal feedback $G^{\ast}$:

\begin{definition}
	Let $G:[0,T]\times L^2(\Lambda)\rightarrow L^2(\Lambda)$ and $n\in \mathbb{N}$. We say that a sequence of subsets $(\mathbb{U}^{n,m})_{m\in \mathbb{N}}\subset U_{\text{ad}}$ satisfies the uniform Lipschitz approximation property with respect to $G$ in dimension $n$, if there exists a sequence of Lipschitz continuous controls $(G^{n,m})_{m\in \mathbb{N}}$ with Lipschitz constants independent of $m$, such that $G^{n,m}\in \mathbb{U}^{n,m}$, and a sequence of radii $(R_m^n)_{m\in \mathbb{N}}$ with $\lim_{m\rightarrow \infty}R_m^n=\infty$, such that
	\begin{align}
		\epsilon_m^n &   :=\sup_{(t,u)\in [0,T]\times \mathcal{B}_{L^2(\Lambda)}(0,R_m^n)}\|P_n(G^n(t,u)-G^{n,m}(t,u))\|_{L^2(\Lambda)}^2 \\ 	\label{RateEpsilon} 
		& \rightarrow 0,
	\end{align}
	as $m\rightarrow \infty$.
\end{definition}
For this type of ansatz spaces, Theorem \ref{main2} provides error estimates for
\begin{align}
	|\inf_{g\in \mathbb{U}^{n,m}}J(g)-\inf_{\mathfrak{g}\in \mathbb{A}}J(\mathfrak{g})|
\end{align}
in terms of $\epsilon^n_m$ and the projection error $\gamma_n$ \eqref{gamma_n}.

\subsection{Universal Approximation}

We start by constructing an ansatz space that satisfies a uniform approximation property with respect to $G^{\ast}$. For $n\in \mathbb{N}$, we consider the function $g^n:[0,T]\times \mathbb{R}^{n}\rightarrow \mathbb{R}^{n}$ given by
\begin{equation}
	g_k^n(t,u_1,\dots,u_n)=\langle G^{\ast}(t,\sum_{j=1}^n u_j e_j),e_k\rangle_{L^2(\Lambda)},\quad k=1,\dots,n.
\end{equation}
Since $G^{\ast}$ is continuous, the functions $(g^n)_{n\in \mathbb{N}}$ are also continuous. In particular, it is possible to approximate these functions by simpler functions that can be treated numerically, e.g., artificial neural networks. In the following we consider for any $n\in \mathbb{N}$ a set $\mathcal{N}^n$ of Lipschitz continuous approximations, such that for all $R>0$ there exists a sequence $(\psi_m^n)_{m\in \mathbb{N}}\in \mathcal{N}^n$ with 
\begin{align}
	\sup_{(t,x)\in [0,T]\times \mathcal{B}_{\mathbb{R}^n}(0,R)}|\psi^n_m(t,x)-g^n(t,x)|^2\rightarrow 0,\quad m\rightarrow \infty,
\end{align}
where 
\begin{align}
	\mathcal{B}_{\mathbb{R}^n}(0,R):=\{x\in \mathbb{R}^n||x|\leq R\}.
\end{align}
For a particular choice of $\mathcal{N}^n$ we refer to our examples in Section \ref{examples}. Then we define the ansatz space 
\begin{multline}
	\mathbb{U}:=\Big \{G(t,\sum_{j=0}^\infty u_je_j)
	=\sum_{k=0}^n\psi(t,\eta^l(u_1,\dots,u_n))_ke_k\\
	\Big |\psi\in \mathcal{N}^n,n,l\in \mathbb{N}\Big\},
\end{multline}
where $\eta^l:\mathbb{R}^n\rightarrow \mathbb{R}^n$ 
\begin{align}
	\eta^l(x)&=\begin{cases}
		x & |x|\leq l\\
		l\frac{x}{|x|} & |x|>l
	\end{cases}
\end{align}
is a smooth cutoff function. It is not difficult to see that for any $n\in \mathbb{N}$ and $R>0$ there exists a sequence $(G^{n,m})_{m\in \mathbb{N}}\subset \mathbb{U}$ that satisfies a linear growth condition of the type 
\begin{align}
	\|G^{n,m}(t,u)\|_{L^2(\Lambda)}&\leq C_n(1+\|u\|_{L^2(\Lambda)}),
\end{align}
such that for any $n\in \mathbb{N}$ 
\begin{align}
	\lim\limits_{m\rightarrow \infty}\sup_{(t,u)\in [0,T]\times \mathcal{B}_{L^2(\Lambda)}(0,R)}\|G^n(t,u)-G^{n,m}(t,u)\|_{L^2(\Lambda)}^2=0,
\end{align}
where 
\begin{align}
	\mathcal{B}_{L^2(\Lambda)}(0,R):=\{u\in L^2(\Lambda)|\|u\|_{L^2(\Lambda)}&\leq R\}.
\end{align}
In particular $\mathbb{U}$ satisfies a uniform approximation property with respect to $G^{\ast}$. We will give a short proof for this in Section \ref{examples} where we consider artificial neural networks, mapping from $\mathbb{R}^{n+1}$ to $\mathbb{R}^n$, as an explicit example for Lipschitz continuous approximations $\mathcal{N}^n$.

Our first main result shows that we can reach the optimal cost of the control problem when we consider ansatz spaces of the above type. A detailed overview of the assumptions for Theorem \ref{main1} can be found in \cite{stannat2023} H1) - H5).

\begin{theorem}\label{main1}
	Assume that $\mathbb{U}\subset U_{\text{ad}}$ satisfies the uniform approximation property with respect to $G^{\ast}$. Then 
	\begin{align}
		\inf_{\mathfrak{g}\in \mathbb{A}}J(\mathfrak{g})=\inf_{G\in \mathbb{U}}J(G).
	\end{align}
\end{theorem}

\subsection{Finite Dimensional Approximation}\label{finiteDimAnsatz}

Our second main result, Theorem \ref{main2}, requires that the controls of our ansatz space take values in the Sobolev space $H^{1}(\Lambda)$. Therefore we strengthen our assumption on the finite dimensional subspaces and assume that 
\begin{align}
	S_n\subset H^{1}(\Lambda).
\end{align}
In order to obtain convergence rates, we also need to specify the rate of convergence of the orthogonal projections. Therefore, we assume that 
\begin{align}\label{gamma_n}
	\|P_nu-u\|_{L^2(\Lambda)}\leq \gamma_n\|u\|_{H^{1}(\Lambda)},
\end{align}
for some $\gamma_n\rightarrow 0$, as $n\rightarrow \infty$. Furthermore, we assume that there exists an optimal control $\mathfrak{g}^{\ast}$ in feedback form with a Lipschitz continuous feedback function $G^{\ast}$. A simple example for such a situation is the linear quadratic case discussed in Example \ref{linearquadratic}.

Since $G^{\ast}$ is Lipschitz continuous, the functions $g^n$ are also Lipschitz continuous. As we will see in the examples in Section \ref{examples}, it is therefore possible to approximate $g^n$ by artificial neural networks that have uniformly bounded Lipschitz constants. For any $n\in \mathbb{N}$, let $(\mathcal{N}^{n,m})_{m\in \mathbb{N}}$ be a sequence of sets of Lipschitz continuous approximations, such that for all $m\in \mathbb{N}$ there exists a $\psi^{n,m}\in \mathcal{N}^{n,m}$ with Lipschitz constant independent of $m$ and
\begin{align}
	\sup_{(t,x)\in [0,T]\times \mathcal{B}_{\mathbb{R}^n}(0,R_m^n)}|\psi^{n,m}(t,x)-g^n(t,x)|^2=: \epsilon^n_m \rightarrow 0,
\end{align}
as $m\rightarrow \infty$, for some sequence $(R_m^n)_{m\in \mathbb{N}}$ of radii with $R_m^n\rightarrow \infty$, as $m\rightarrow \infty$. Then we define the sequence of ansatz spaces of dimension $n\in \mathbb{N}$ by 
\begin{multline}
	\mathbb{U}^{n,m}:=\Big \{G(t,\sum_{j=0}^\infty u_je_j) =\sum_{k=0}^n\psi(t,u_1,\dots,u_n)_ke_k\\
	\Big |\psi\in \mathcal{N}^{n,m}\Big\}.
\end{multline}
It is not difficult to observe that for any $m\in \mathbb{N}$, there exists a Lipschitz continuous control $G^{n,m}\in \mathbb{U}^{n,m}$ with Lipschitz constant independent of $m$, such that 
\begin{align}
	\sup_{(t,u)\in [0,T]\times \mathcal{B}_{L^2(\Lambda)}(0,R^n_m)}\|P_n(G^{n,m}(t,u)-G^{n}(t,u))\|_{L^2(\Lambda)}^2\rightarrow 0.
\end{align}
In particular, for any $n\in \mathbb{N}$, the sequence $(\mathbb{U}^{n,m})_{m\in \mathbb{N}}$ satisfies the Lipschitz approximation property with respect to $G^{\ast}$ in dimension $n$.

Our second main result provides convergence rates for ansatz spaces of the above type. A detailed overview of the assumptions for Theorem \ref{main2} can be found in \cite{stannat2023} H1) - H5) and S1) - S3).

\begin{theorem}\label{main2}
	Let $n\in \mathbb{N}$. Let $(\mathbb{U}^{n,m})_{m\in \mathbb{N}}$ be a sequence of subsets of $U_{\text{ad}}$ that satisfies the Lipschitz approximation property with respect to $G^{\ast}$ in dimension $n$. Then it holds 
	\begin{multline}
		\inf_{G\in \mathbb{U}^{n,m}}J(G)-\inf_{\mathfrak{g}\in \mathbb{A}}J(\mathfrak{g})\\
		\leq C\left (1+\sqrt{\mathbb{E} \left [ \int_{0}^{T}\|\mathfrak{g}^{\ast}_t\|_{H^1(\Lambda)}^2\mathrm{d}t \right ]}\right )\gamma_n+C_n\sqrt{\epsilon^n_m+\frac{1}{R^n_m}},
	\end{multline}
	for some universal constant $C$ independent of $n$ and $m$ and some constant $C_n$ which is independent of $m$. Here $\epsilon_m^n$ is given by \eqref{RateEpsilon}. 
\end{theorem}

\section{Gradient Descent Algorithm}\label{gradientdescent}

In this section, we describe our gradient descent algorithm. Theorem \ref{main1} enables us to consider the approximating optimal control problem of minimizing optimal costs on a finite dimensional ansatz space $\mathbb{U}^{n,m}$. We assume that 
\begin{equation}
	\mathbb{U}^{n,m} = \left \{ \Phi(\cdot,\cdot,\alpha) \Big | \alpha \in \mathbb{R}^{d_m} \right \}
\end{equation}
for a parametrization $\Phi:[0,T]\times L^2(\Lambda)\times \mathbb{R}^{d_m}\to L^2(\Lambda)$. Replacing
\begin{equation}
	G(t,u) = \Phi(t,u,\alpha)
\end{equation}
for some given $\alpha\in\mathbb{R}^{d_m}$ leads to the state equation
\begin{equation}\label{statealpha}
	\begin{cases}
		\mathrm{d}u^{\alpha}_t = [\Delta u^{\alpha}_t + \mathcal{F}(u^{\alpha}_t)+\Phi(t,u^{\alpha}_t,\alpha) ] \mathrm{d}t +\sigma \mathrm{d}W_t\\
		u^{\alpha}_0=u\in L^2(\Lambda),
	\end{cases}
\end{equation}
and the cost functional $J:\mathbb{R}^{d_m}\to\mathbb{R}$,
\begin{multline}
	J(\alpha) = \mathbb{E}\bigg [ \int_0^T \int_{\Lambda} l(t,x,u^{\alpha}_t(x)) \mathrm{d}x + \frac{\nu}{2} \| \Phi(t,u^{\alpha}_t,\alpha) \|_{L^2(\Lambda)}^2 \mathrm{d}t\\
	+ \int_{\Lambda} m(x,u^{\alpha}_T(x)) \mathrm{d}x \bigg ].
\end{multline}
Using this parametrization, we obtain under suitable regularity assumptions on $\Phi$ the following representation of the gradient of the cost functional:
\begin{multline}\label{gradient}
	\nabla J(\alpha)\\
	= \mathbb{E} \left [ \int_0^T \nu \Phi^{\ast}_{\alpha}(t,u^{\alpha}_t,\alpha) \Phi(t,u^{\alpha}_t,\alpha) + \Phi^{\ast}_{\alpha}(t,u^{\alpha}_t,\alpha) p_t \mathrm{d}t \right ],
\end{multline}
where $p:[0,T]\times \Lambda\times \Omega \to \mathbb{R}$ is the solution of the so-called adjoint equation
\begin{equation}\label{adjoint}
	\begin{cases}
		\mathrm{d}p_t = - [ (\Delta + \mathcal{F}^{\prime}(u^{\alpha}_t) + \Phi^{\ast}_u(t,u^{\alpha}_t,\alpha))p_t + \mathcal{L}^{\prime}(t,u^{\alpha}_t)\\
		\qquad\qquad\qquad\qquad\qquad\quad + \nu \Phi_u^{\ast}(t,u^{\alpha}_t,\alpha)\Phi(t,u^{\alpha}_t,\alpha)] \mathrm{d}t\\
		p_T=\mathcal{M}^{\prime}(u^{\alpha}_T).
	\end{cases}
\end{equation}
Here
\begin{equation}
	\mathcal{L}(t,u) := \int_{\Lambda}l(t,x,u(x))\mathrm{d}x, \quad u\in L^2(\Lambda),t\in [0,T]
\end{equation}
and
\begin{equation}
	\mathcal{M}(u) := \int_{\Lambda}m(x,u(x))\mathrm{d}x, \quad u\in L^2(\Lambda).
\end{equation}
Furthermore, $\mathcal{F}^{\prime}$, $\mathcal{L}^{\prime}$ and $\mathcal{M}^{\prime}$, denote the associated Fr\'echet derivatives with respect to $u$. Finally, $\Phi_{\alpha}(t,u,\alpha):\mathbb{R}^{d_m}\to L^2(\Lambda)$ (resp. $\Phi_u(t,u,\alpha):L^2(\Lambda)\to L^2(\Lambda)$) denotes the derivative of $\Phi$ with respect to $\alpha$ (resp. $u$), and $\Phi_{\alpha}^{\ast}(t,u,\alpha)$ (resp. $\Phi_u^{\ast}(t,u,\alpha)$) denotes the adjoint of $\Phi_{\alpha}(t,u,\alpha)$ (resp. $\Phi_u(t,u,\alpha)$). For a derivation of the gradient, see the Appendix. 

Note that equation \eqref{adjoint} is a linear backward PDE with random coefficients which are given by the state equation \eqref{statealpha}. For details concerning the numerical implementation, we refer to our software made available on GitHub, see \cite{vogler2023}. 

\begin{example}[Artificial Neural Network]
	In the case of an artificial neural network with activator function $\theta$, the controls in the ansatz space $\mathbb{U}^{n,m}$ can be parametrized as
	\begin{equation}
		\Phi(t,u,\alpha) = \sum_{i=1}^n \left ( B \theta \left (A \begin{pmatrix} t\\ \pi_n u \end{pmatrix} + a \right ) \right )_i e_i,
	\end{equation}
	where $\alpha=(A,a,B)$ consists of $A\in \mathbb{R}^{k\times (n+1)}$, $a\in \mathbb{R}^{k}$, and $B\in \mathbb{R}^{n \times k}$, for respective dimensions $n=n_m$ and $k=k_m$, and $\pi_n u=(\langle u,e_1\rangle,\dots,\langle u,e_n \rangle )^{\top}$.
\end{example}

Based on this representation, we implement the following algorithm:
\begin{algorithm}\label{gradientdescentalgorithm}
	Fix an initial control $\alpha_0$, a stopping criterion $\rho>0$, and a step size $s>0$.
	\begin{enumerate}
		\item\label{stepone} Solve the state equation \eqref{statealpha} for one realization of the noise.
		\item\label{steptwo} Solve the adjoint equation \eqref{adjoint} with the data given by the sample calculated in Step \ref{stepone}.
		\item Repeat Step \ref{stepone} and Step \ref{steptwo} to approximate the gradient \eqref{gradient} using Monte Carlo approximation.
		\item Compute new control via $\alpha_{n+1}=\alpha_n-s \nabla J(\alpha_n)$.
		\item Stop if $\|\nabla J(\alpha_n)\|<\rho$.
	\end{enumerate}
\end{algorithm}

\section{Artificial Neural Networks as Ansatz Spaces}\label{examples}

In this section, we discuss explicit examples for ansatz spaces $\mathbb{U}\subset U_{\text{ad}}$ that satisfy the assumptions of our main results, Theorem \ref{main1} and Theorem \ref{main2}, and are suitable for the numerical implementation of our gradient descent algorithm. In the first part of this section, we focus on universal spaces and Theorem \ref{main1}. The second part is devoted to ansatz spaces of fixed size and Theorem \ref{main2} with corresponding convergence rates. In the whole section, we consider the finite dimensional subspaces 
\begin{equation}
	S_n:=\text{span}\left \{ 1,\sqrt{2}\cos(k \pi\,\cdot\,)\middle |\,k=1,\dots,n \right \} \subset H^{1}(0,1), 
\end{equation}
with orthonormal basis in $L^2(\Lambda)$
\begin{align}
	e_0=1, \quad e_k=\sqrt{2}\cos(k \pi \,\cdot\,),\quad k=1,\dots,n.
\end{align}
In particular we have for any $u\in H^{1}(0,1)$
\begin{align}
	&\|P_nu-u\|_{L^2(0,1)}^2\\
    &=\sum_{k=n+1}^\infty |\langle u,\sqrt{2}\cos(k \pi \,\cdot\,)\rangle_{L^2(0,1)}|^2\\
	&\leq\frac{1}{\pi^2 n^2} \sum_{k=n+1}^\infty |\langle (-\Delta)^{1/2} u,\sqrt{2}\cos(k \pi \,\cdot\,)\rangle_{L^2(0,1)}|^2\\
	&\leq \frac{1}{\pi^2 n^2} \|u\|_{H^{1}(0,1)}^2=:\gamma_n\|u\|_{H^{1}(0,1)}^2,
\end{align}
where we used that the $k$-th eigenvalue of the Neumann Laplace operator with respect to the eigenfunction $e_k$ is given by $-\pi^2 k^2$.

\subsection{Neural Network Approximation}
\label{sec:NeuralNetworkApproximation} 

Regarding Theorem \ref{main1}, we will show in our first example that it is indeed sufficient to consider the type of ansatz space constructed in the first part of Section \ref{Ansatz}, using one-layer artificial neural networks for the approximating sets, to get arbitrarily close to the optimal cost. More precisely, we show that the set 
\begin{multline}
	\mathbb{U}:=\bigg \{G(t,\sum_{i=1}^\infty u_i e_i)(x):=\sum_{i=1}^n\psi(t,\eta^l(u_1,\dots,u_n))_ie_i(x)\\
	\bigg |\psi\in \mathcal{N}^n,n,l\in \mathbb{N} \bigg \}
\end{multline}
satisfies the uniform approximation property with respect to $G^{\ast}$, where 
\begin{align}
	\mathcal{N}^n:=\bigcup_{k=1}^\infty \mathcal{N}^n_k
\end{align}
and 
\begin{multline}
	\mathcal{N}_k^n:=\big \{\psi(x)=B \theta (A x+a)\\
	\big |\,A\in \mathbb{R}^{k\times (n+1)},a\in \mathbb{R}^k,B\in \mathbb{R}^{(n+1)\times k}\big \}
\end{multline}
denotes the set of all one-layer artificial neural networks from $\mathbb{R}^n$ to $\mathbb{R}^n$ with $k$ neurons, for a given non-polynomial, Lipschitz continuous activator function $\theta$.

To this end, we recall the following classical universal approximation result by \cite{Pinkus1999}:

\begin{theorem}\label{thm:main3}
	Let $\theta\in \mathcal{C}(\mathbb{R})$, then we define for $u=(u_1,\dots,u_d)\in \mathbb{R}^d$
	\begin{equation}
		\theta(u)^i:=\theta(u_i).
	\end{equation}
	If $\theta$ is not polynomial, then for any $n,m\in \mathbb{N}$, compact set $K\subset \mathbb{R}^n$, $h\in \mathcal{C}(K,\mathbb{R}^m)$ and $\epsilon>0$, there exists $k\in \mathbb{N},A\in \mathbb{R}^{k\times n},a\in \mathbb{R}^k,B\in \mathbb{R}^{m\times k}$, such that 
	\begin{align}
		\sup_{u\in K}|h(u)-\psi(u)|<\epsilon,
	\end{align}
	where $\psi$ is the one-layer artificial neural network 
	\begin{align}
		\psi(u):=B \theta (A u+a).
	\end{align}
\end{theorem}

Let $n\in \mathbb{N}$ and recall from Section \ref{Ansatz} the finitely based approximation of $G^{\ast}$
\begin{align}
	g^n&:[0,T]\times \mathbb{R}^n\rightarrow \mathbb{R}^n\\
	g^n_i(t,u)&:=\langle G^{\ast}(t,\sum_{j=1}^{n}u_je_j),e_i\rangle,\quad i=1,\dots, n.
\end{align}
Since $g^n$ is continuous, there exists for any $m\in \mathbb{N}$ a one-layer artificial neural network $\psi^{n,m}:[0,T]\times \mathbb{R}^n\rightarrow \mathbb{R}^n$, such that 
\begin{align}
	\sup_{(t,u)\in [0,T]\times \mathcal{B}_{\mathbb{R}^n}(0,m)}|g^n(t,u)-\psi^{n,m}(t,u)|<\frac{1}{m}.
\end{align}
Since $G^{\ast}$ satisfies a linear growth condition, i.e.,
\begin{align}
	\|G^{\ast}(t,u)\|_{L^2(\Lambda)}\leq C(1+\|u\|_{L^2(\Lambda)}),
\end{align}
for some $C>0$, we have 
\begin{align}
	|g^n(t,u)|&=\|P_n G^{\ast}(t,\sum_{i=1}^{n}u_i e_i)\|_{L^2(\Lambda)}\\
	&\leq \|G^{\ast}(t,\sum_{i=1}^{n}u_ie_i)\|_{L^2(\Lambda)}\\
	&\leq C(1+|u|),
\end{align}
where $u=(u_1,\dots,u_n)$. Hence, if we consider the continuous function $\eta^m:\mathbb{R}^n\rightarrow \mathbb{R}^n$ 
\begin{align}
	\eta^m(x)&=\begin{cases}
		x & |x|\leq m\\
		m\frac{x}{|x|} & |x|>m
	\end{cases}
\end{align}
and define 
\begin{align}
	\tilde{\psi}^{n,m}(t,u):=\psi^{n,m}(t,\eta(u)),
\end{align}
then clearly $\tilde{\psi}^{n,m}=\psi^{n,m}$ on $[0,T]\times \mathcal{B}_{\mathbb{R}^n}(0,m)$ and for any $(t,u)\in [0,T]\times \mathcal{B}_{\mathbb{R}^n}(0,m)$
\begin{align}
	|\psi^{n,m}(t,u)|&\leq |g^n(t,u)|+|\psi^{n,m}(t,u)-g^n(t,u)|\\
	&\leq C(1+|u|)+1.
\end{align}
Furthermore, on $[0,T]\times  \mathcal{B}_{\mathbb{R}^n}(0,m)^c$ we have 
\begin{align}
	|\tilde{\psi}^{n,m}(t,u)|&=|\psi^{n,m}(t,m\frac{u}{|u|})|\\
	&\leq |g^n(t,m\frac{u}{|u|})|+|\psi^{n,m}(t,m\frac{u}{|u|})-g^n(t,m\frac{u}{|u|})|\\
	&\leq C(1+m)+1\\
	&\leq C(1+|u|)+1.
\end{align}
Therefore $\tilde{\psi}^{n,m}$ satisfies a linear growth condition with some constant $C>0$ independent of $m$. Now, we define
\begin{align}
	G^{n,m}(t,u):=\sum_{i=1}^{n}\tilde{\psi}^{n,m}(t,\langle u,e_1\rangle_{L^2(\Lambda)},\ldots,\langle u,e_n\rangle_{L^2(\Lambda)}))_ie_i.
\end{align}
One can easily check that all the assumptions of Theorem \ref{main2} are satisfied for $\mathbb{U}$. Indeed, due to the Lipschitz continuity of the elements in $\mathcal{N}_k^n$, any control in $\mathbb{U}$ is indeed admissible, i.e., $\mathbb{U}\subset U_{\text{ad}}$. Furthermore, $(G^{n,m})_{m\in \mathbb{N}}$ is a sequence in $\mathbb{U}$ that satisfies a linear growth condition with some constant $C>0$ independent of $m$ and approximates any finitely based $G^n$. Indeed, for any $R>0$ and any $\epsilon>0$ there exists an $M\in \mathbb{N}$, such that $ \mathcal{B}_{L^2(\Lambda)}(0,R)\subset  \mathcal{B}_{L^2(\Lambda)}(0,m)$ and $\frac{1}{m}<\epsilon$ for every $m\geq M$. Therefore, we have for any $m\geq M$
\begin{align}
	&\sup_{(t,u)\in [0,T]\times \mathcal{B}_{L^2(\Lambda)}(0,R)}\|G^{n,m}(t,u)-G^n(t,u)\|^2\\
	&\leq\sup_{(t,u)\in [0,T]\times \mathcal{B}_{L^2(\Lambda)}(0,m)}\|G^{n,m}(t,u)-G^n(t,u)\|^2\\
	&\leq \sup_{(t,u)\in [0,T]\times \mathcal{B}_{\mathbb{R}^n}(0,m)}|\tilde{\psi}^{n,m}(t,u)-g^n(t,u)|^2\\
	&\leq \sup_{(t,u)\in [0,T]\times \mathcal{B}_{\mathbb{R}^n}(0,m)}|\psi^{n,m}(t,u)-g^n(t,u)|^2\\
	&<\frac{1}{m}<\epsilon.
\end{align}

In the case of bounded controls, for example if $\mathcal{U}=\mathcal{B}_{L^2(\Lambda)}(0,R)$, we could consider the ansatz space
\begin{multline}
	\mathbb{U}:=\bigg \{G(t,\sum_{i=1}^\infty u_i e_i)(x):=\sum_{i=1}^n\psi(t,\eta^l(u_1,\dots,u_n))_ie_i(x)\\
	\bigg |\|G_n(t,u)\|_{L^2(\Lambda)}\leq R\text{ where }\psi\in \mathcal{N}^n,n,l\in \mathbb{N} \bigg \}
\end{multline}
and the sequence 
\begin{multline}
	G^{n,m}(t,u)\\
	:=\!\left (\!1-\frac{1}{m(R+1)}\!\right )\!\!\sum_{i=1}^{n}\!\tilde{\psi}^{n,m}(t,\langle u,e_1\rangle_{L^2(\Lambda)},\dots,\langle u,e_n\rangle_{L^2(\Lambda)}))_ie_i.
\end{multline}

\subsection{Convergence Rates}\label{section:convrates}

Next, we provide explicit convergence rates for the sequence of approximating spaces
\begin{multline}
	\mathbb{U}^{n,m}:=\Big\{G(t,\sum_{i=1}^\infty u_ie_i)(x)=\sum_{i=1}^n\psi(t,(u_1,\dots,u_n))_ie_i(x)\\
	\Big |\psi\in \mathcal{N}_m^n\Big\}
\end{multline}
using Theorem \ref{main2}. We will mainly follow the ideas of \cite{carmona2022}. In the following we consider a $2\pi$-periodic activator function $\theta:\mathbb{R}\rightarrow \mathbb{R}$ that satisfies $\hat{\theta}:=\int_{-\pi}^{\pi}\theta(x)e^{-\mathrm{i}x}\mathrm{d}x\not=0$. Recall the following result from \cite{carmona2022}:

\begin{proposition}\label{neuralapprox}
	Let $K,L_1,L_2,R>0$ and $h:[0,T]\times \mathcal{B}_{\mathbb{R}^{n+1}}(0,R)\rightarrow \mathbb{R}^n$ be Lipschitz continuous in $(t,x)$ with $\|h\|_{\mathcal{C}^0([0,T]\times \mathcal{B}_{\mathbb{R}^{n+1}}(0,R))}\leq K$ and Lipschitz constant bounded by $L_1$. Furthermore we assume that $h$ is twice differentiable in $x$ and $\partial_{x_i}h$ is Lipschitz continuous with Lipschitz constant bounded by $L_2$. Then there exists a constant $C>0$ depending only on the above constants, on $n,T$ and on the activator function $\theta$ through $\hat{\theta}$, $\|\theta'\|_{\mathcal{C}^0}$, $\|\theta''\|_{\mathcal{C}^0}$ and $\|\theta'''\|_{\mathcal{C}^0}$, and there exists a constant $m_0$ depending only on $n$ with the following property. For every $m_{\text{in}}>m_0$, there exists a one-hidden layer artificial neural network $\psi_h\in \mathcal{N}_{m_{\text{in}}}^n$ such that 
	\begin{align}
		\|h-\psi_h\|_{\mathcal{C}^0([0,T]\times \mathcal{B}_{\mathbb{R}^{n}}(0,R);\mathbb{R}^n)}\leq C(1+R)m_{\text{in}}^{-1/(2(n+1))}
	\end{align}
	and such that the Lipschitz constants of $\psi_h,\partial_x\psi_h$ are at most $C(1+Rm_{\text{in}}^{-1/(2(n+1))})$.
\end{proposition}
If we now set $R^n_m:=m^{1/(3(n+1))}$, then by Proposition \ref{neuralapprox} there exists for any sufficiently large $m\in \mathbb{N}$ a network $\psi^{n,m}\in \mathcal{N}_{m}^n$, such that for
\begin{equation}
	G^{n,m}(t,\sum_{i=1}^\infty u_i e_i)(x)=\sum_{i=1}^n\psi^{n,m}(t,(u_1,\dots,u_n))_ie_i(x),
\end{equation}
it holds
\begin{multline}
	\sup_{(t,x)\in [0,T]\times \mathcal{B}_{L^2(\Lambda)}(0,R^n(m))}\|P_n(G^{n,m}(t,x)-G^{n}(t,x))\|_{L^2(\Lambda)}^2\\
	\leq m^{-1/(3(n+1))}\rightarrow 0.
\end{multline}
Furthermore $G^{n,m}\in \mathbb{U}^{n,m}$ is Lipschitz continuous with some Lipschitz constant independent of $n$. Therefore, applying Theorem \ref{main2} yields 
\begin{multline}
	\inf_{G\in \mathbb{U}^{n,m}}J(G)-\inf_{\mathfrak{g}\in \mathbb{A}}J(\mathfrak{g})\\
	\leq C\left (1+\sqrt{\mathbb{E}\left[\int_{0}^{T}\|\mathfrak{g}^{\ast}_t\|_{H^1(\Lambda)}^2\mathrm{d}t\right]}\right )\frac{1}{n^2\pi^2}+C_n\sqrt{m^{-1/(3(n+1))}}.
\end{multline}

\section{Finite Element Discretization of the Control Problem}\label{numericaldiscretization}

In order to implement a numerical algorithm, the control problem needs to be discretized. In this section, we introduce the finite element discretized version of our control problem and provide a bound for the error resulting from the discretization.

\begin{problem}[Finite Element Control Problem]\label{FESCP}
	Let $n\in \mathbb{N}$. Minimize 
	\begin{align}
		J^n(G)&:=\mathbb{E}\left[\int_{0}^{T}\int_{\Lambda} l(t,x,u^{G,n}_t(x)) \mathrm{d}x + \frac{\nu}{2} \| P_nG(t,u^{G,n}_t) \|^2_{L^2(\Lambda)} \mathrm{d}t\right] \\
		&\quad + \mathbb{E}\left[\int_{\Lambda}m(x,u^{G,n}_T(x))\mathrm{d}x \right]
	\end{align}
	over the set $U_{\text{ad}}$ subject to the discretized SPDE
	\begin{equation}\label{FEspde}
		\begin{cases}
			\mathrm{d}u^{G,n}_t=[ -\Delta_n u^{{G},n}_t+P_n \mathcal{F}(u^{{G},n}_t)+P_n{G}(t,u^{{G},n}_t)]\mathrm{d}t+ \sigma  P_n \mathrm{d}W_t\\
			u^{{G},n}_0=P_n u\in S_n,
		\end{cases}
	\end{equation}
	where $\Delta_n x_n$ is defined as the unique element in $S_n$ with 
	\begin{align}
		\langle \nabla x_n,\nabla y_n \rangle_{L^2(\Lambda)}=\langle \Delta_n x_n,y_n\rangle_{L^2(\Lambda)}\quad \forall y_n\in S_n.
	\end{align}
\end{problem}

\begin{proposition}
	\label{prop:numericaldiscretization}
	Under the same assumptions as in Theorem \ref{main2}, it holds 
	\begin{multline}
		\inf_{G\in \mathbb{U}^{n,m}}J^n(G)-\inf_{\mathfrak{g}\in \mathbb{A}}J(\mathfrak{g})\\
		\leq C\left (1+\sqrt{\mathbb{E} \left [ \int_{0}^{T}\|\mathfrak{g}^{\ast}_t\|_{H^1(\Lambda)}^2\mathrm{d}t \right ]}\right )\gamma_n+C_n\sqrt{\epsilon^n_m+\frac{1}{R_m^n}},
	\end{multline}
	for some constant $C_n>0$ that only depends on $n\in \mathbb{N}$ and on the Lipschitz constant of $G^{\ast}$.
\end{proposition}

\begin{remark}
	Under additional assumptions, in particular convexity of the Hamiltonian, one can show that for any $\mathbb{U}\subset U_{\text{ad}}$, it holds 
	\begin{align}
		\inf_{\mathfrak{g}\in \mathbb{A}}J(\mathfrak{g})-\inf_{G\in \mathbb{U}}J^n(G)\leq C\gamma_n,
	\end{align}
	for some constant $C>0$ which is independent of $n$.
\end{remark}

\section{Simulations}\label{simulations}

In this whole section, we consider 
equation \eqref{spde1} on some interval 
$\Lambda = (0,L)$ of length $L$. 
For the approximation we consider 
the Galerkin finite dimensional subspace
\begin{equation}
	S_n:=\text{span}\left\{ 1, \sqrt{\frac{2}{L}}\cos\left( \frac{k}{L} \pi\,\cdot\, \right) \middle |\,  k=1,\dots,n\right \}\subset H^{1}(0,L), 
\end{equation}
with orthonormal basis
\begin{equation}
	e_0=\sqrt{\frac{1}{L}},\quad e_k =\sqrt{\frac{2}{L}}\cos\left( \frac{k}{L} \pi \,\cdot\right) ,\quad k=1,\dots,n,
\end{equation}
or the finite element subspace
\begin{equation}
	\overline{S}_n:=\text{span}\left \{ \overline{e}_k\middle |\,k=1,\dots,n \right \} \subset H^{1}(0,L), 
\end{equation}
with basis
\begin{equation}
	\overline{e}_k(x)=\begin{cases}
		n\left (x-\frac{k-1}{n}L\right ) &\text{if } x\in \left [\frac{k-1}{n}L,\frac{k}{n}L \right ]\\
		n\left (\frac{k+1}{n}L -x\right ) & \text{if } x\in \left [\frac{k}{n}L,\frac{k+1}{n}L \right ] \\
		0 & \text{otherwise} 
	\end{cases}
\end{equation}
for $k = 1, \dots , n$. In the following we will always consider $n=400$ for the simulations and solve the finite-element discretized version of equation \eqref{statealpha} and \eqref{adjoint} using the semi-implicit Euler-Maruyama method with time discretization $\Delta t=0.05$.

\subsection{Heat Equation}

In this subsection, we consider the controlled stochastic heat equation in order to validate our algorithm by comparing with the optimal feedback control obtained from the associated Riccati equation, see Example \ref{linearquadratic}. The controlled state equation is given by
\begin{equation}\label{A:stateequation}
	\begin{cases}
		\mathrm{d}u^{\mathfrak{g}}_t = [\Delta u^{\mathfrak{g}}_t + \mathfrak{g}_t]  \mathrm{d}t + 0.05\mathrm{d}W_t,\quad t\in [0,20]\\
		u^{\mathfrak{g}}_0=u\in L^2(0,20),
	\end{cases}
\end{equation}
with Neumann boundary conditions and $u = \mathbf{1}_{[20/3,40/3]}$. We consider the problem of steering the solution of the stochastic heat equation into the constant zero profile. To this end, we introduce the cost functional
\begin{equation}
	J(\mathfrak{g}) = \frac{1}{2} \mathbb{E} \left [ \int_0^{20} \| u^{\mathfrak{g}}_t\|_{L^2(0,20)}^2 + \| \mathfrak{g}_t \|_{L^2(0,20)}^2 \mathrm{d}t \right ].
\end{equation} 
Note that the second term is a regularization, which is necessary in linear quadratic control theory. We approximate the Riccati equation \eqref{Ricc} numerically based on $400$ Fourier coefficients to obtain the approximated optimal feedback control 
\begin{equation}
	\mathfrak{g}^{\text{Ric}}_t = P^n(t) u^{\mathfrak{g}^{\text{Ric}}}_t
\end{equation}
and use this approximation as a benchmark for our gradient descent algorithm.

For our approximation we consider the ansatz space $\mathbb{U}^{n,m}$ constructed in Section \ref{section:convrates} using artificial neural networks with $m=50$ neurons.

After about $40,000$ iterations of our stochastic gradient descent algorithm, we end up with an approximated cost of $J\approx 4.58$, and approximated $L^2$-distance of our approximation $\mathfrak{g}^{\text{approx}}_t=G^{\text{approx}}(t,u_t^{G^{\text{approx}}})$ to the optimal control $\mathfrak{g}^{\text{Ric}}$ given by
\begin{align}
	\mathbb{E}\left[\int_{0}^{20} \|\mathfrak{g}_t^{\text{Ric}}-\mathfrak{g}^{\text{approx}}_t\|_{L^2(\Lambda)}^2\mathrm{d}t\right]\approx 0.02.
\end{align}

\begin{figure}[H]
	\centering
	\includegraphics[width=1.0\linewidth]{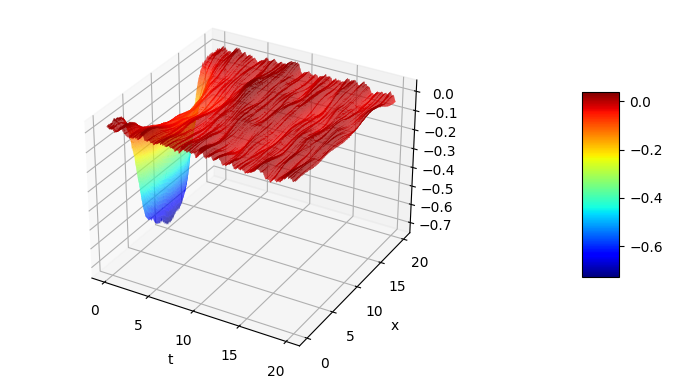}
	\caption{Sample of approximated optimal control $\mathfrak{g}^{\text{approx}}$}\label{A:feedbackcontrol}
\end{figure}

\begin{figure}[H]
	\centering
	\includegraphics[width=1.0\linewidth]{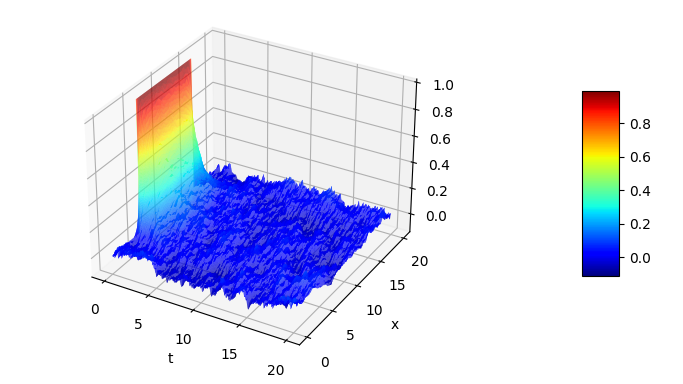}
	\caption{Sample of controlled solution $u^{\mathfrak{g}^{\text{approx}}}$}\label{A:controlledstate}
\end{figure}

\begin{figure}[H]
	\centering
	\includegraphics[width=1.0\linewidth]{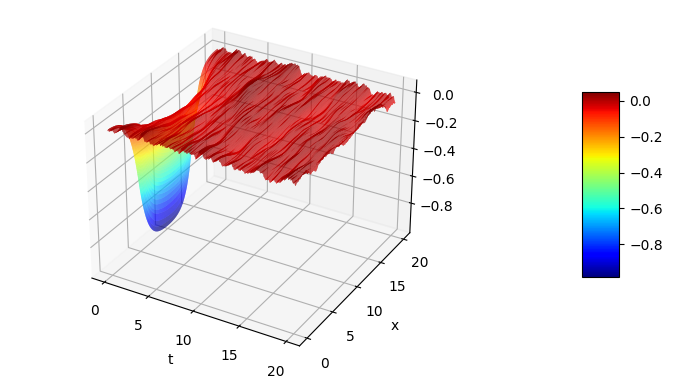}
	\caption{Sample of optimal control $\mathfrak{g}^{\text{Ric}}$}\label{A:riccaticontrol}
\end{figure}

Below, we display our simulation results. Figure \ref{A:uncontrolledstate} displays one realization of the uncontrolled stochastic heat equation. Figure \ref{A:feedbackcontrol} displays our neural network approximation of the optimal feedback control, and Figure \ref{A:controlledstate} shows its impact when applied to the system \eqref{A:stateequation}. Figure \ref{A:riccaticontrol} displays the optimal control obtained using the Riccati equation, and shows that our approximated feedback control is indeed qualitatively close.

\begin{figure}[H]
	\centering
	\includegraphics[width=1.0\linewidth]{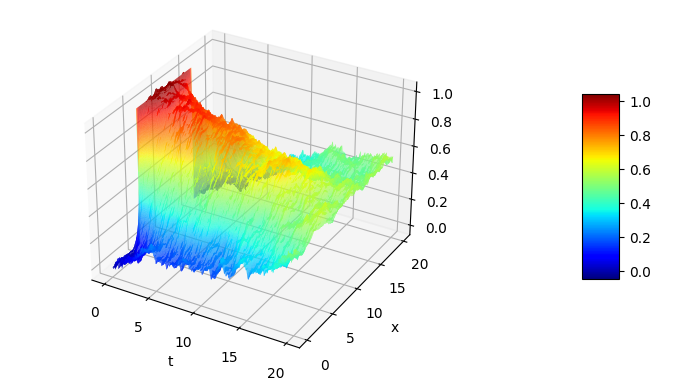}
	\caption{Sample of uncontrolled solution}\label{A:uncontrolledstate}
\end{figure}

\subsection{$L^2$-Feedback Control of the Nagumo Equation}\label{NagumoFeedback}

\begin{figure}[H]
	\centering
	\includegraphics[width=1\linewidth]{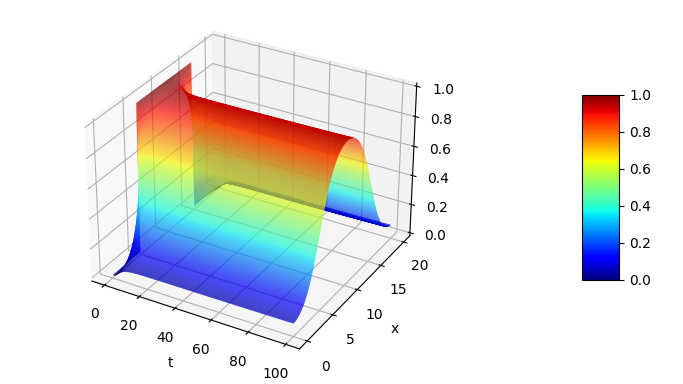}
	\caption{Reference profile $u^0$}\label{B:referenceprofile}
\end{figure}

In this example, we apply our algorithm to the controlled stochastic Nagumo equation
\begin{equation}\label{Nagumo}
	\begin{cases}
		\mathrm{d}u^G_t \! = \! \left [ \Delta u^G_t - u^G_t (u^G_t-\frac12 )(u^G_t-1)+G(t,u^G_t) \right ]\! \mathrm{d}t + \! 0.05\mathrm{d}W_t\\
		u^G_0 = u\in L^2(0,20),
	\end{cases}
\end{equation}
with Neumann boundary conditions and $u= \mathbf{1}_{[5,15]}$. We consider the problem of stabilizing a bump profile given by the solution to the uncontrolled deterministic Nagumo equation, i.e., equation \eqref{Nagumo} with $G\equiv 0$ and without noise, see Figure \ref{B:referenceprofile}. To this end, we introduce the cost functional
\begin{align}
	J(G) &= \mathbb{E} \left [ \int_0^{100} \| u^{G}_t- u^0_t\|_{L^2(0,20)}^2 + \frac{1}{2}  \| G(t,u^G_t) \|_{L^2(0,20)}^2 \mathrm{d}t \right ]\\
	&\quad + \mathbb{E} \left [ \| u^G_{100}- u^0_{100}\|_{L^2(0,20)}^2 \right ].\label{costfunctionalB}
\end{align} 
For the approximation of the optimal control, we use the ansatz space 
\begin{multline}\label{ansatzNN}
	\mathbb{U}^{n,k_1+k_2}:=\bigg\{G(t,\sum_{i=1}^\infty u_ie_i)(x)=\sum_{i=1}^n\psi(t,(u_1,\dots,u_n))_ie_i(x)\\
	\bigg |\psi\in \mathcal{N}_{k_1,k_2}^n\bigg\},
\end{multline}
for $k=k_1=k_2=100$, where 
\begin{multline}
	\mathcal{N}_{k,k}^n:=\bigg \{\psi(t,u)=C\theta\left (B \theta \left (A \begin{pmatrix}t\\ u
	\end{pmatrix}+a\right )+b \right )\\
	\bigg |\,A\in \mathbb{R}^{k\times (n+1)},B\in \mathbb{R}^{k\times k},C\in \mathbb{R}^{(n+1)\times k},a,b\in \mathbb{R}^{k} \bigg \}
\end{multline}
denotes the set of all two-layer neural networks with $k$ neurons in the first and second layer, and ReLU activator function $\theta$. In this case, our optimal control achieves an approximated cost of $J\approx 8.1$.

\begin{figure}[H]
	\centering
	\includegraphics[width=1\linewidth]{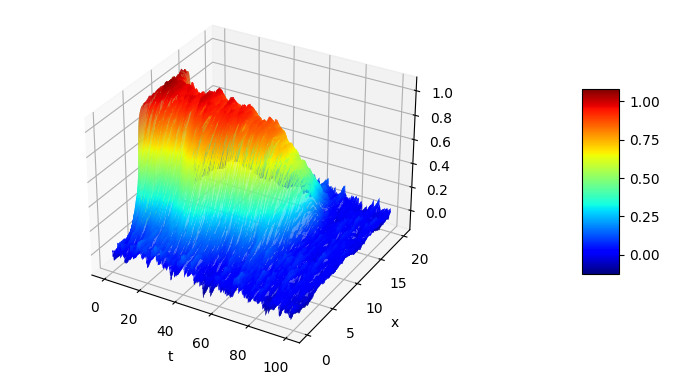}
	\caption{Sample of uncontrolled state}\label{B:sampletwo}
\end{figure}

\begin{figure}[H]
	\centering
	\includegraphics[width=1\linewidth]{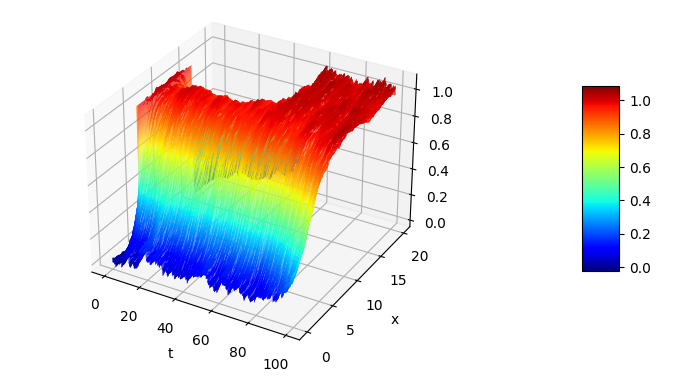}
	\caption{Sample of uncontrolled state}\label{B:sampleone}
\end{figure}

\begin{figure}[H]
	\centering
	\includegraphics[width=1\linewidth]{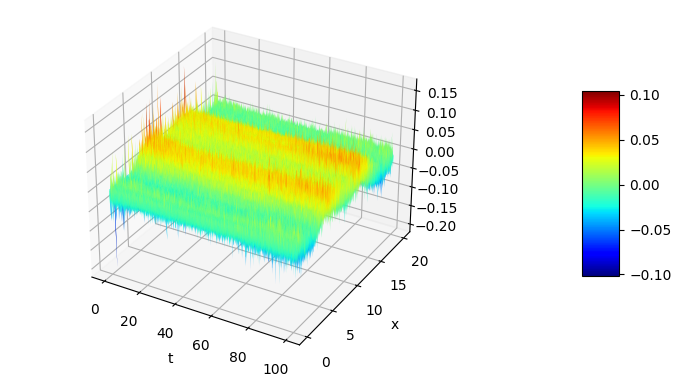}
	\caption{Sample of approximated optimal control}\label{B:optimalcontrol}
\end{figure}

\begin{figure}[H]
	\centering
	\includegraphics[width=1\linewidth]{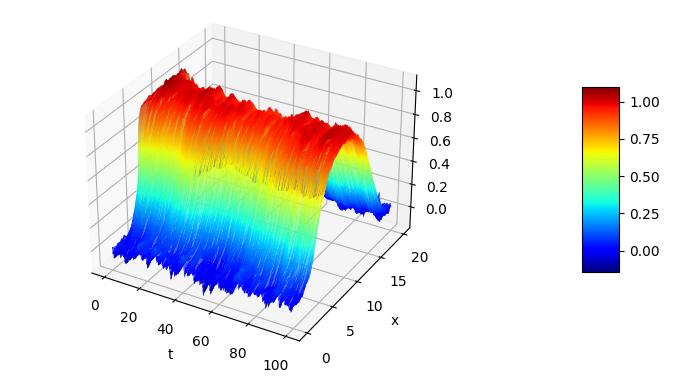}
	\caption{Sample of controlled state}\label{B:controlledstate}
\end{figure}

Figures \ref{B:sampletwo} and \ref{B:sampleone} display two realizations of the uncontrolled stochastic Nagumo equation \eqref{Nagumo}, i.e., $G\equiv 0$. Without control, the bump is unstable and the noise pushes the solution to one of the stable steady states $u\equiv 0$ or $u\equiv 1$. Figure \ref{B:optimalcontrol} displays one realization of the approximated optimal control and shows that the control mostly acts on the interface, but also reacts to the noise in the system. Figure \ref{B:controlledstate} illustrates the impact of the approximated optimal control when applied to the system \eqref{Nagumo}; it shows that our feedback control indeed stabilizes the bump.

\subsection{Nemytskii Feedback Control of the Nagumo Equation}
In our final example, we consider again the controlled stochastic Nagumo equation \eqref{Nagumo} with the same cost functional \eqref{costfunctionalB}. However, now we only consider feedback controls of Nemytskii-type, i.e., the control $G$ is a Nemytskii operator 
\begin{align}
	\mathcal{G}(t,u)(x)&=g(t,x,u(x)) 
\end{align}
for some function $g:[0,100]\times [0,20] \times \mathbb{R}\to\mathbb{R}$. This means that the control at point $x\in\Lambda$ depends on $u(x)$, the value of the solution at point $x$, but not on the whole function $u\in L^2(\Lambda)$ as in the previous example. This restriction to feedback controls of Nemytskii-type significantly reduces the computational complexity and therefore leads to a more efficient approximation.

\begin{figure}[H]
	\centering
	\includegraphics[width=1\linewidth]{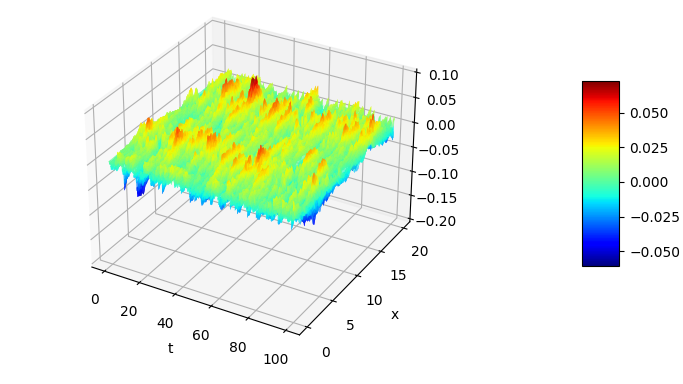}
	\caption{Sample of approximated optimal control}\label{C:optimalcontrol}
\end{figure}

\begin{figure}[H]
	\centering
	\includegraphics[width=1\linewidth]{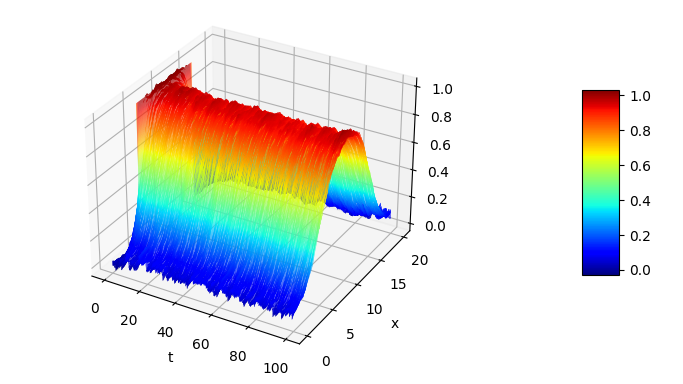}
	\caption{Sample of controlled state}\label{C:controlledstate}
\end{figure}

\begin{figure}[H]
	\centering
	\includegraphics[width=1\linewidth]{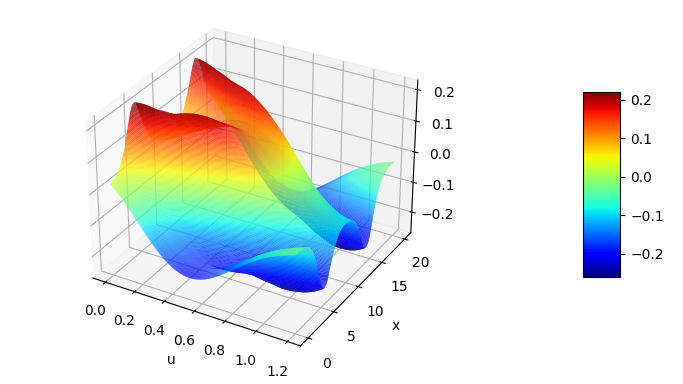}
	\caption{Approximation of feedback function $(x,u)\mapsto g(t,x,u)$ at time $t=40$}\label{C:feedbackfunction}
\end{figure}

For the approximation of the optimal control, we consider the following ansatz space of Gaussian radial basis function neural networks with $m=40$ neurons:
\begin{align}
	&\mathbb{U}^{n,m}\\
	&:=\bigg \{G(t,u)(x) = \sum_{i=1}^{n} \sum_{j=1}^m \sum_{k=1}^r \alpha_{ijk} \mathbf{1}_{[t_{k-1}, t_k]}(t) e^{-\kappa |u(x)-\overline{u}_j|^2}  \overline{e}_i(x)\\
	&\qquad\qquad\qquad\qquad\qquad\qquad\qquad\qquad \bigg |\alpha_{ijk}\in \mathbb{R}, \overline{u}_j\in \mathbb{R}\bigg \},
\end{align}
where $\{0=t_0<t_1<\ldots<t_r=100\}$ and $\kappa = 6$. In this case, our optimal control achieves an approximated cost of $J\approx 1.25$. Similar to the approximated control of Section \ref{NagumoFeedback}, Figure \ref{C:optimalcontrol} shows that the control again mostly acts on the interface and also reacts to the noise in the system (compare with Figure \ref{B:optimalcontrol}). Figure \ref{C:controlledstate} shows that the Nemytskii feedback control achieves a better result than the $L^2$-feedback control from Subsection \ref{NagumoFeedback} (compare with Figure \ref{B:controlledstate}). Figure \ref{C:feedbackfunction} displays the approximation of the feedback function $(x,u) \mapsto g(t,x,u)$ at time $t=40$. Observe that the feedback function is only trained for profiles arising in simulations. In particular, the feedback function is not trained in the middle of the bump ($x\approx 10$) for $u\approx 0$.

\section{Conclusion and Outlook}\label{conclusion}

In this work we presented a direct method to approximate optimal feedback controls for stochastic reaction diffusion equations using finitely based approximations combined with artificial neural networks. We were able to validate our method numerically in the linear quadratic case. Furthermore our numerical results indicate that the method is able to approximate efficient controls even in more complex situations.
On the theoretical side we provided approximation results under rather weak assumptions on the control problem using one-layer artificial neural networks. To the best of our knowledge there are no results concerning the approximation by neural networks, that would help improving our overall approximation rate. However we expect better approximation rates using more involved types of neural network architecture, such as deep neural networks, which can be investigated in future work. We also expect that our method can be applied to mean-field control problems.

%\section*{Supplemental Material}
%TBD

\begin{acknowledgments}
	This work was funded by the Deutsche Forschungsgemeinschaft (DFG) via grant
	CRC 910, ``Control of self-organizing nonlinear systems: Theoretical methods and concepts of application,'' project A10, ``Control of stochastic mean-field equations with applications to brain networks.''
\end{acknowledgments}

\section*{Author Declarations}
\subsection*{Conflict of Interest}
The authors have no conflicts to disclose.
\subsection*{Author Contributions}
Wilhelm Stannat: Mathematical theory (equal). Alexander Vogler: Mathematical theory (equal); Software. Lukas Wessels: Mathematical theory (equal).

\section*{Data Availability}
Software and simulation output freely available with publication on GitHub at \url{https://github.com/AVoglerTu/SFB910Feedback}, see \cite{vogler2023}.

\appendix*

\section{Derivation of the Gradient of the Cost Functional}

As discussed in Section V, we now consider the following approximating optimal control problem: Minimize
\begin{multline}
	J(\alpha) = \mathbb{E}\Big [ \int_0^T \int_{\Lambda} l(t,x,u^{\alpha}_t(x)) \mathrm{d}x + \frac{\nu}{2} \| \Phi(t,u^{\alpha}_t,\alpha) \|_{L^2(\Lambda)}^2 \mathrm{d}t\\
	+ \int_{\Lambda} m(x,u^{\alpha}_T(x)) \mathrm{d}x \Big ].
\end{multline}
subject to
\begin{equation}\label{statealphaappendix}
	\begin{cases}
		\mathrm{d}u^{\alpha}_t = [\Delta u^{\alpha}_t + \mathcal{F}(u^{\alpha}_t)+\Phi(t,u^{\alpha}_t,\alpha) ] \mathrm{d}t +\sigma \mathrm{d}W_t\\
		u^{\alpha}_0=u\in L^2(\Lambda).
	\end{cases}
\end{equation}
In the following derivation of the gradient of the cost functional, we assume that $l$ and $m$ are differentiable with respect to $u$, $\mathcal{F}$ is differentiable, and $\Phi$ is differentiable with respect to $u$ and $\alpha$. 
The derivation of the gradient of the cost functional follows along the same lines as the derivation in \cite{stannat20212}, however, in the present setting, the control enters the state equation \eqref{statealphaappendix} in a nonlinear fashion, which requires slight modifications of the arguments. We begin by differentiating the cost functional in some direction $\beta\in \mathbb{R}^{d_m}$ which yields
\begin{align}\label{directionalderivative}
	&\frac{\partial J(\alpha)}{\partial \beta}\\
	&= \mathbb{E} \Big [ \int_0^T \langle \mathcal{L}^{\prime} (t,u^{\alpha}_t), y^{\beta}_t \rangle + \nu \langle \Phi(t,u^{\alpha}_t,\alpha), \Phi_{\alpha}(t,u^{\alpha}_t,\alpha) \beta \rangle\\
	&\quad +\nu \langle \Phi(t,u^{\alpha}_t,\alpha),\Phi_u(t,u^{\alpha}_t,\alpha) y^{\beta}_t\rangle \mathrm{d}t +\langle \mathcal{M}^{\prime}(u^{\alpha}_T), y^{\beta}_T \rangle \Big ],
\end{align}
where $y^{\beta}$ is the solution of the linearized state equation
\begin{equation}
	\begin{cases}
		\mathrm{d}y^{\beta}_t = [ \Delta y^{\beta}_t + \mathcal{F}^{\prime} (u^{\alpha}_t) y^{\beta}_t + \Phi_u(t,u^{\alpha},\alpha) y^{\beta}_t \\
		\qquad\qquad\qquad\qquad\qquad\qquad\quad +\Phi_{\alpha}(t,u^{\alpha}_t,\alpha) \beta ] \mathrm{d}t\\
		y^{\beta}_0=0.
	\end{cases}
\end{equation}
Next, we introduce the adjoint state $p$ as the solution of the adjoint equation
\begin{equation}
	\begin{cases}
		\mathrm{d}p_t = - [ (\Delta + \mathcal{F}^{\prime}(u^{\alpha}_t) + \Phi^{\ast}_u(t,u^{\alpha}_t,\alpha))p_t + \mathcal{L}^{\prime}(t,u^{\alpha}_t)\\
		\qquad\qquad\qquad\qquad\qquad\quad + \nu \Phi_u^{\ast}(t,u^{\alpha}_t,\alpha)\Phi(t,u^{\alpha}_t,\alpha)] \mathrm{d}t\\
		p_T=\mathcal{M}^{\prime}(u^{\alpha}_T).
	\end{cases}
\end{equation}
A straightforward computation leads to the following adjoint state property:
\begin{align}
	\mathrm{d}\langle y^{\beta}_t,p_t \rangle &= \langle y^{\beta}_t,\mathrm{d}p_t \rangle + \langle p_t,\mathrm{d}y^{\beta}_t \rangle\\
	&=\big [ -\nu \langle \Phi^{\ast}_u(t,u^{\alpha}_t,\alpha)\Phi(t,u^{\alpha}_t,\alpha),y^{\beta}_t \rangle - \langle \mathcal{L}^{\prime}(t,u^{\alpha}_t),y^{\beta}_t \rangle\\
	&\qquad+ \langle \Phi_{\alpha}(t,u^{\alpha}_t,\alpha)\beta, p_t \rangle \big ]\mathrm{d}t.
\end{align}
Integrating over $[0,T]$ and taking expectations, we obtain from equation \eqref{directionalderivative} the desired representation of the gradient of the cost functional:
\begin{multline}
	\nabla J(\alpha)\\
	= \mathbb{E} \left [ \int_0^T \nu \Phi^{\ast}_{\alpha}(t,u^{\alpha}_t,\alpha) \Phi(t,u^{\alpha}_t,\alpha) + \Phi^{\ast}_{\alpha}(t,u^{\alpha}_t,\alpha) p_t \mathrm{d}t \right ].
\end{multline}

\section*{References}

\nocite{*}
\bibliography{StannatVoglerWessels}% Produces the bibliography via BibTeX.

\begin{thebibliography}{40}%
\makeatletter
\providecommand \@ifxundefined [1]{%
 \@ifx{#1\undefined}
}%
\providecommand \@ifnum [1]{%
 \ifnum #1\expandafter \@firstoftwo
 \else \expandafter \@secondoftwo
 \fi
}%
\providecommand \@ifx [1]{%
 \ifx #1\expandafter \@firstoftwo
 \else \expandafter \@secondoftwo
 \fi
}%
\providecommand \natexlab [1]{#1}%
\providecommand \enquote  [1]{``#1''}%
\providecommand \bibnamefont  [1]{#1}%
\providecommand \bibfnamefont [1]{#1}%
\providecommand \citenamefont [1]{#1}%
\providecommand \href@noop [0]{\@secondoftwo}%
\providecommand \href [0]{\begingroup \@sanitize@url \@href}%
\providecommand \@href[1]{\@@startlink{#1}\@@href}%
\providecommand \@@href[1]{\endgroup#1\@@endlink}%
\providecommand \@sanitize@url [0]{\catcode `\\12\catcode `\$12\catcode
  `\&12\catcode `\#12\catcode `\^12\catcode `\_12\catcode `\%12\relax}%
\providecommand \@@startlink[1]{}%
\providecommand \@@endlink[0]{}%
\providecommand \url  [0]{\begingroup\@sanitize@url \@url }%
\providecommand \@url [1]{\endgroup\@href {#1}{\urlprefix }}%
\providecommand \urlprefix  [0]{URL }%
\providecommand \Eprint [0]{\href }%
\providecommand \doibase [0]{http://dx.doi.org/}%
\providecommand \selectlanguage [0]{\@gobble}%
\providecommand \bibinfo  [0]{\@secondoftwo}%
\providecommand \bibfield  [0]{\@secondoftwo}%
\providecommand \translation [1]{[#1]}%
\providecommand \BibitemOpen [0]{}%
\providecommand \bibitemStop [0]{}%
\providecommand \bibitemNoStop [0]{.\EOS\space}%
\providecommand \EOS [0]{\spacefactor3000\relax}%
\providecommand \BibitemShut  [1]{\csname bibitem#1\endcsname}%
\let\auto@bib@innerbib\@empty
%</preamble>
\bibitem [{\citenamefont {Showalter}(1997)}]{showalter1997}%
  \BibitemOpen
  \bibfield  {author} {\bibinfo {author} {\bibfnamefont {R.~E.}\ \bibnamefont
  {Showalter}},\ }\href@noop {} {\emph {\bibinfo {title} {Monotone operators in
  {B}anach space and nonlinear partial differential equations}}}\ (\bibinfo
  {publisher} {American Mathematical Society},\ \bibinfo {address} {Providence,
  RI},\ \bibinfo {year} {1997})\BibitemShut {NoStop}%
\bibitem [{\citenamefont {Liu}\ and\ \citenamefont
  {R\"ockner}(2015)}]{liu2015}%
  \BibitemOpen
  \bibfield  {author} {\bibinfo {author} {\bibfnamefont {W.}~\bibnamefont
  {Liu}}\ and\ \bibinfo {author} {\bibfnamefont {M.}~\bibnamefont
  {R\"ockner}},\ }\href@noop {} {\emph {\bibinfo {title} {Stochastic Partial
  Differential Equations: An Introduction}}}\ (\bibinfo  {publisher}
  {Springer},\ \bibinfo {year} {2015})\BibitemShut {NoStop}%
\bibitem [{\citenamefont {Hinze}\ \emph {et~al.}(2009)\citenamefont {Hinze},
  \citenamefont {Pinnau}, \citenamefont {Ulbrich},\ and\ \citenamefont
  {Ulbrich}}]{hinze2009}%
  \BibitemOpen
  \bibfield  {author} {\bibinfo {author} {\bibfnamefont {M.}~\bibnamefont
  {Hinze}}, \bibinfo {author} {\bibfnamefont {R.}~\bibnamefont {Pinnau}},
  \bibinfo {author} {\bibfnamefont {M.}~\bibnamefont {Ulbrich}}, \ and\
  \bibinfo {author} {\bibfnamefont {S.}~\bibnamefont {Ulbrich}},\ }\href@noop
  {} {\emph {\bibinfo {title} {Optimization with {PDE} constraints}}}\
  (\bibinfo  {publisher} {Springer},\ \bibinfo {year} {2009})\BibitemShut
  {NoStop}%
\bibitem [{\citenamefont {Li}\ and\ \citenamefont {Yong}(1995)}]{li1995}%
  \BibitemOpen
  \bibfield  {author} {\bibinfo {author} {\bibfnamefont {X.}~\bibnamefont
  {Li}}\ and\ \bibinfo {author} {\bibfnamefont {J.}~\bibnamefont {Yong}},\
  }\href@noop {} {\emph {\bibinfo {title} {Optimal Control Theory for Infinite
  Dimensional Systems}}}\ (\bibinfo  {publisher} {Birkh\"auser},\ \bibinfo
  {address} {Boston, MA},\ \bibinfo {year} {1995})\BibitemShut {NoStop}%
\bibitem [{\citenamefont {Tr\"oltzsch}(2010)}]{troeltzsch2010}%
  \BibitemOpen
  \bibfield  {author} {\bibinfo {author} {\bibfnamefont {F.}~\bibnamefont
  {Tr\"oltzsch}},\ }\href@noop {} {\emph {\bibinfo {title} {Optimal Control of
  Partial Differential Equations}}}\ (\bibinfo  {publisher} {American
  Mathematical Society},\ \bibinfo {address} {Providence, RI},\ \bibinfo {year}
  {2010})\BibitemShut {NoStop}%
\bibitem [{\citenamefont {Buchholz}\ \emph
  {et~al.}(2013{\natexlab{a}})\citenamefont {Buchholz}, \citenamefont {Engel},
  \citenamefont {Kammann},\ and\ \citenamefont {Tr\"oltzsch}}]{buchholz2013}%
  \BibitemOpen
  \bibfield  {author} {\bibinfo {author} {\bibfnamefont {R.}~\bibnamefont
  {Buchholz}}, \bibinfo {author} {\bibfnamefont {H.}~\bibnamefont {Engel}},
  \bibinfo {author} {\bibfnamefont {E.}~\bibnamefont {Kammann}}, \ and\
  \bibinfo {author} {\bibfnamefont {F.}~\bibnamefont {Tr\"oltzsch}},\
  }\bibfield  {title} {\enquote {\bibinfo {title} {{O}n the optimal control of
  the {S}chl\"{o}gl-model},}\ }\href {\doibase 10.1007/s10589-013-9550-y}
  {\bibfield  {journal} {\bibinfo  {journal} {Comput. Optim. Appl.}\ }\textbf
  {\bibinfo {volume} {56}},\ \bibinfo {pages} {153--185} (\bibinfo {year}
  {2013}{\natexlab{a}})}\BibitemShut {NoStop}%
\bibitem [{\citenamefont {Buchholz}\ \emph
  {et~al.}(2013{\natexlab{b}})\citenamefont {Buchholz}, \citenamefont {Engel},
  \citenamefont {Kammann},\ and\ \citenamefont {Tr\"oltzsch}}]{buchholz20132}%
  \BibitemOpen
  \bibfield  {author} {\bibinfo {author} {\bibfnamefont {R.}~\bibnamefont
  {Buchholz}}, \bibinfo {author} {\bibfnamefont {H.}~\bibnamefont {Engel}},
  \bibinfo {author} {\bibfnamefont {E.}~\bibnamefont {Kammann}}, \ and\
  \bibinfo {author} {\bibfnamefont {F.}~\bibnamefont {Tr\"oltzsch}},\
  }\bibfield  {title} {\enquote {\bibinfo {title} {{E}rratum to: {O}n the
  optimal control of the {S}chl\"{o}gl-model},}\ }\href {\doibase
  10.1007/s10589-013-9570-7} {\bibfield  {journal} {\bibinfo  {journal}
  {Comput. Optim. Appl.}\ }\textbf {\bibinfo {volume} {56}},\ \bibinfo {pages}
  {187--188} (\bibinfo {year} {2013}{\natexlab{b}})}\BibitemShut {NoStop}%
\bibitem [{\citenamefont {Ryll}(2017)}]{ryll2016}%
  \BibitemOpen
  \bibfield  {author} {\bibinfo {author} {\bibfnamefont {C.}~\bibnamefont
  {Ryll}},\ }\emph {\bibinfo {title} {Optimal control of patterns in some
  reaction-diffusion-systems}},\ \href@noop {} {\bibinfo {type} {Doctoral
  thesis}},\ \bibinfo  {school} {Technische Universität Berlin}, \bibinfo
  {address} {Berlin, Germany} (\bibinfo {year} {2017})\BibitemShut {NoStop}%
\bibitem [{\citenamefont {Ryll}\ \emph {et~al.}(2016)\citenamefont {Ryll},
  \citenamefont {L\"ober}, \citenamefont {Martens}, \citenamefont {Engel},\
  and\ \citenamefont {Tr\"oltzsch}}]{ryll2017}%
  \BibitemOpen
  \bibfield  {author} {\bibinfo {author} {\bibfnamefont {C.}~\bibnamefont
  {Ryll}}, \bibinfo {author} {\bibfnamefont {J.}~\bibnamefont {L\"ober}},
  \bibinfo {author} {\bibfnamefont {S.}~\bibnamefont {Martens}}, \bibinfo
  {author} {\bibfnamefont {H.}~\bibnamefont {Engel}}, \ and\ \bibinfo {author}
  {\bibfnamefont {F.}~\bibnamefont {Tr\"oltzsch}},\ }\bibfield  {title}
  {\enquote {\bibinfo {title} {Analytical, optimal, and sparse optimal control
  of traveling wave solutions to reaction-diffusion systems},}\ }in\ \href@noop
  {} {\emph {\bibinfo {booktitle} {Control of self-organizing nonlinear
  systems}}},\ \bibinfo {editor} {edited by\ \bibinfo {editor} {\bibfnamefont
  {E.}~\bibnamefont {Sch\"oll}}, \bibinfo {editor} {\bibfnamefont {S.~H.~L.}\
  \bibnamefont {Klapp}}, \ and\ \bibinfo {editor} {\bibfnamefont
  {P.}~\bibnamefont {H\"ovel}}}\ (\bibinfo  {publisher} {Springer},\ \bibinfo
  {year} {2016})\ pp.\ \bibinfo {pages} {189--210}\BibitemShut {NoStop}%
\bibitem [{\citenamefont {Cordoni}\ and\ \citenamefont
  {Di~Persio}(2018)}]{cordoni2018}%
  \BibitemOpen
  \bibfield  {author} {\bibinfo {author} {\bibfnamefont {F.}~\bibnamefont
  {Cordoni}}\ and\ \bibinfo {author} {\bibfnamefont {L.}~\bibnamefont
  {Di~Persio}},\ }\bibfield  {title} {\enquote {\bibinfo {title} {Optimal
  control for the stochastic {F}itz{H}ugh-{N}agumo model with recovery
  variable},}\ }\href {\doibase 10.3934/eect.2018027} {\bibfield  {journal}
  {\bibinfo  {journal} {Evol. Eq. Control Theory}\ }\textbf {\bibinfo {volume}
  {7}},\ \bibinfo {pages} {571--585} (\bibinfo {year} {2018})}\BibitemShut
  {NoStop}%
\bibitem [{\citenamefont {Du}\ and\ \citenamefont {Meng}(2013)}]{du2013}%
  \BibitemOpen
  \bibfield  {author} {\bibinfo {author} {\bibfnamefont {K.}~\bibnamefont
  {Du}}\ and\ \bibinfo {author} {\bibfnamefont {Q.}~\bibnamefont {Meng}},\
  }\bibfield  {title} {\enquote {\bibinfo {title} {A maximum principle for
  optimal control of stochastic evolution equations},}\ }\href {\doibase
  10.1137/120882433} {\bibfield  {journal} {\bibinfo  {journal} {SIAM J.
  Control Optim.}\ }\textbf {\bibinfo {volume} {51}},\ \bibinfo {pages}
  {4343--4362} (\bibinfo {year} {2013})}\BibitemShut {NoStop}%
\bibitem [{\citenamefont {Fabbri}, \citenamefont {Gozzi},\ and\ \citenamefont
  {{\'{S}}wi{\k{e}}ch}(2017)}]{fabbri2017}%
  \BibitemOpen
  \bibfield  {author} {\bibinfo {author} {\bibfnamefont {G.}~\bibnamefont
  {Fabbri}}, \bibinfo {author} {\bibfnamefont {F.}~\bibnamefont {Gozzi}}, \
  and\ \bibinfo {author} {\bibfnamefont {A.}~\bibnamefont
  {{\'{S}}wi{\k{e}}ch}},\ }\href@noop {} {\emph {\bibinfo {title} {Stochastic
  Optimal Control in Infinite Dimension}}}\ (\bibinfo  {publisher} {Springer},\
  \bibinfo {year} {2017})\BibitemShut {NoStop}%
\bibitem [{\citenamefont {Frankowska}\ and\ \citenamefont
  {Zhang}(2020)}]{frankowska2020}%
  \BibitemOpen
  \bibfield  {author} {\bibinfo {author} {\bibfnamefont {H.}~\bibnamefont
  {Frankowska}}\ and\ \bibinfo {author} {\bibfnamefont {X.}~\bibnamefont
  {Zhang}},\ }\bibfield  {title} {\enquote {\bibinfo {title} {Necessary
  conditions for stochastic optimal control problems in infinite dimensions},}\
  }\href {\doibase 10.1016/j.spa.2019.11.010} {\bibfield  {journal} {\bibinfo
  {journal} {Stochastic Process. Appl.}\ }\textbf {\bibinfo {volume} {130}},\
  \bibinfo {pages} {4081--4103} (\bibinfo {year} {2020})}\BibitemShut {NoStop}%
\bibitem [{\citenamefont {Fuhrman}, \citenamefont {Hu},\ and\ \citenamefont
  {Tessitore}(2013)}]{fuhrman2013}%
  \BibitemOpen
  \bibfield  {author} {\bibinfo {author} {\bibfnamefont {M.}~\bibnamefont
  {Fuhrman}}, \bibinfo {author} {\bibfnamefont {Y.}~\bibnamefont {Hu}}, \ and\
  \bibinfo {author} {\bibfnamefont {G.}~\bibnamefont {Tessitore}},\ }\bibfield
  {title} {\enquote {\bibinfo {title} {Stochastic maximum principle for optimal
  control of {S}{P}{D}{E}s},}\ }\href {\doibase 10.1007/s00245-013-9203-7}
  {\bibfield  {journal} {\bibinfo  {journal} {Appl. Math. Optim.}\ }\textbf
  {\bibinfo {volume} {68}},\ \bibinfo {pages} {181--217} (\bibinfo {year}
  {2013})}\BibitemShut {NoStop}%
\bibitem [{\citenamefont {Fuhrman}, \citenamefont {Hu},\ and\ \citenamefont
  {Tessitore}(2018)}]{fuhrman2018}%
  \BibitemOpen
  \bibfield  {author} {\bibinfo {author} {\bibfnamefont {M.}~\bibnamefont
  {Fuhrman}}, \bibinfo {author} {\bibfnamefont {Y.}~\bibnamefont {Hu}}, \ and\
  \bibinfo {author} {\bibfnamefont {G.}~\bibnamefont {Tessitore}},\ }\bibfield
  {title} {\enquote {\bibinfo {title} {Stochastic maximum principle for optimal
  control of partial differential equations driven by white noise},}\ }\href
  {\doibase 10.1007/s40072-017-0108-3} {\bibfield  {journal} {\bibinfo
  {journal} {Stoch. Partial Differ. Equ. Anal. Comput.}\ }\textbf {\bibinfo
  {volume} {6}},\ \bibinfo {pages} {255--285} (\bibinfo {year}
  {2018})}\BibitemShut {NoStop}%
\bibitem [{\citenamefont {Fuhrman}\ and\ \citenamefont
  {Orrieri}(2016)}]{fuhrman2016}%
  \BibitemOpen
  \bibfield  {author} {\bibinfo {author} {\bibfnamefont {M.}~\bibnamefont
  {Fuhrman}}\ and\ \bibinfo {author} {\bibfnamefont {C.}~\bibnamefont
  {Orrieri}},\ }\bibfield  {title} {\enquote {\bibinfo {title} {Stochastic
  maximum principle for optimal control of a class of nonlinear {S}{P}{D}{E}s
  with dissipative drift},}\ }\href {\doibase 10.1137/15M1012888} {\bibfield
  {journal} {\bibinfo  {journal} {SIAM J. Control Optim.}\ }\textbf {\bibinfo
  {volume} {54}},\ \bibinfo {pages} {341--371} (\bibinfo {year}
  {2016})}\BibitemShut {NoStop}%
\bibitem [{\citenamefont {L{\"u}}\ and\ \citenamefont {Zhang}(2014)}]{lue2014}%
  \BibitemOpen
  \bibfield  {author} {\bibinfo {author} {\bibfnamefont {Q.}~\bibnamefont
  {L{\"u}}}\ and\ \bibinfo {author} {\bibfnamefont {X.}~\bibnamefont {Zhang}},\
  }\href@noop {} {\emph {\bibinfo {title} {General {P}ontryagin-Type Stochastic
  Maximum Principle and Backward Stochastic Evolution Equations in Infinite
  Dimensions}}}\ (\bibinfo  {publisher} {Springer},\ \bibinfo {year}
  {2014})\BibitemShut {NoStop}%
\bibitem [{\citenamefont {L{\"u}}\ and\ \citenamefont {Zhang}(2015)}]{lue2015}%
  \BibitemOpen
  \bibfield  {author} {\bibinfo {author} {\bibfnamefont {Q.}~\bibnamefont
  {L{\"u}}}\ and\ \bibinfo {author} {\bibfnamefont {X.}~\bibnamefont {Zhang}},\
  }\bibfield  {title} {\enquote {\bibinfo {title} {Transposition method for
  backward stochastic evolution equations revisited, and its application},}\
  }\href {\doibase 10.3934/mcrf.2015.5.529} {\bibfield  {journal} {\bibinfo
  {journal} {Math. Control Relat. Fields}\ }\textbf {\bibinfo {volume} {5}},\
  \bibinfo {pages} {529--555} (\bibinfo {year} {2015})}\BibitemShut {NoStop}%
\bibitem [{\citenamefont {L{\"u}}\ and\ \citenamefont {Zhang}(2018)}]{lue2018}%
  \BibitemOpen
  \bibfield  {author} {\bibinfo {author} {\bibfnamefont {Q.}~\bibnamefont
  {L{\"u}}}\ and\ \bibinfo {author} {\bibfnamefont {X.}~\bibnamefont {Zhang}},\
  }\bibfield  {title} {\enquote {\bibinfo {title} {Operator-valued backward
  stochastic {L}yapunov equations in infinite dimensions, and its
  application},}\ }\href {\doibase 10.3934/mcrf.2018014} {\bibfield  {journal}
  {\bibinfo  {journal} {Math. Control Relat. Fields}\ }\textbf {\bibinfo
  {volume} {8}},\ \bibinfo {pages} {337--381} (\bibinfo {year}
  {2018})}\BibitemShut {NoStop}%
\bibitem [{\citenamefont {Stannat}\ and\ \citenamefont
  {Wessels}(2021{\natexlab{a}})}]{stannat2021}%
  \BibitemOpen
  \bibfield  {author} {\bibinfo {author} {\bibfnamefont {W.}~\bibnamefont
  {Stannat}}\ and\ \bibinfo {author} {\bibfnamefont {L.}~\bibnamefont
  {Wessels}},\ }\bibfield  {title} {\enquote {\bibinfo {title} {Peng's maximum
  principle for stochastic partial differential equations},}\ }\href@noop {}
  {\bibfield  {journal} {\bibinfo  {journal} {SIAM J. Control Optim.}\ }\textbf
  {\bibinfo {volume} {59}},\ \bibinfo {pages} {3552--3573} (\bibinfo {year}
  {2021}{\natexlab{a}})}\BibitemShut {NoStop}%
\bibitem [{\citenamefont {Stannat}\ and\ \citenamefont
  {Wessels}(2022)}]{stannat2022}%
  \BibitemOpen
  \bibfield  {author} {\bibinfo {author} {\bibfnamefont {W.}~\bibnamefont
  {Stannat}}\ and\ \bibinfo {author} {\bibfnamefont {L.}~\bibnamefont
  {Wessels}},\ }\href@noop {} {\enquote {\bibinfo {title} {Necessary and
  sufficient conditions for optimal control of semilinear stochastic partial
  differential equations},}\ }\bibinfo {howpublished} {preprint,
  https://arxiv.org/abs/2112.09639} (\bibinfo {year} {2022})\BibitemShut
  {NoStop}%
\bibitem [{\citenamefont {Wessels}(2022)}]{wessels2022}%
  \BibitemOpen
  \bibfield  {author} {\bibinfo {author} {\bibfnamefont {L.}~\bibnamefont
  {Wessels}},\ }\emph {\bibinfo {title} {Optimal control of stochastic
  reaction-diffusion equations}},\ \href@noop {} {\bibinfo {type} {Doctoral
  thesis}},\ \bibinfo  {school} {Technische Universität Berlin}, \bibinfo
  {address} {Berlin, Germany} (\bibinfo {year} {2022})\BibitemShut {NoStop}%
\bibitem [{\citenamefont {Beck}\ \emph {et~al.}(2021)\citenamefont {Beck},
  \citenamefont {Becker}, \citenamefont {Cheridito}, \citenamefont {Jentzen},\
  and\ \citenamefont {Neufeld}}]{beck2021}%
  \BibitemOpen
  \bibfield  {author} {\bibinfo {author} {\bibfnamefont {C.}~\bibnamefont
  {Beck}}, \bibinfo {author} {\bibfnamefont {S.}~\bibnamefont {Becker}},
  \bibinfo {author} {\bibfnamefont {P.}~\bibnamefont {Cheridito}}, \bibinfo
  {author} {\bibfnamefont {A.}~\bibnamefont {Jentzen}}, \ and\ \bibinfo
  {author} {\bibfnamefont {A.}~\bibnamefont {Neufeld}},\ }\bibfield  {title}
  {\enquote {\bibinfo {title} {Deep splitting method for parabolic {PDE}s},}\
  }\href@noop {} {\bibfield  {journal} {\bibinfo  {journal} {SIAM J. Sci.
  Comput.}\ }\textbf {\bibinfo {volume} {43}},\ \bibinfo {pages} {A3135--A3154}
  (\bibinfo {year} {2021})}\BibitemShut {NoStop}%
\bibitem [{\citenamefont {Beck}, \citenamefont {E},\ and\ \citenamefont
  {Jentzen}(2019)}]{beck2019}%
  \BibitemOpen
  \bibfield  {author} {\bibinfo {author} {\bibfnamefont {C.}~\bibnamefont
  {Beck}}, \bibinfo {author} {\bibfnamefont {W.}~\bibnamefont {E}}, \ and\
  \bibinfo {author} {\bibfnamefont {A.}~\bibnamefont {Jentzen}},\ }\bibfield
  {title} {\enquote {\bibinfo {title} {Machine learning approximation
  algorithms for high-dimensional fully nonlinear partial differential
  equations and second-order backward stochastic differential equations},}\
  }\href@noop {} {\bibfield  {journal} {\bibinfo  {journal} {J. Nonlinear
  Sci.}\ }\textbf {\bibinfo {volume} {29}},\ \bibinfo {pages} {1563--1619}
  (\bibinfo {year} {2019})}\BibitemShut {NoStop}%
\bibitem [{\citenamefont {Dolgov}, \citenamefont {Kalise},\ and\ \citenamefont
  {Kunisch}(2021)}]{dolgov2021}%
  \BibitemOpen
  \bibfield  {author} {\bibinfo {author} {\bibfnamefont {S.}~\bibnamefont
  {Dolgov}}, \bibinfo {author} {\bibfnamefont {D.}~\bibnamefont {Kalise}}, \
  and\ \bibinfo {author} {\bibfnamefont {K.~K.}\ \bibnamefont {Kunisch}},\
  }\bibfield  {title} {\enquote {\bibinfo {title} {Tensor decomposition methods
  for high-dimensional {H}amilton-{J}acobi-{B}ellman equations},}\ }\href@noop
  {} {\bibfield  {journal} {\bibinfo  {journal} {SIAM J. Sci. Comput.}\
  }\textbf {\bibinfo {volume} {43}},\ \bibinfo {pages} {A1625--A1650} (\bibinfo
  {year} {2021})}\BibitemShut {NoStop}%
\bibitem [{\citenamefont {Dunst}\ \emph {et~al.}(2019)\citenamefont {Dunst},
  \citenamefont {Majee}, \citenamefont {Prohl},\ and\ \citenamefont
  {Vallet}}]{dunst2019}%
  \BibitemOpen
  \bibfield  {author} {\bibinfo {author} {\bibfnamefont {T.}~\bibnamefont
  {Dunst}}, \bibinfo {author} {\bibfnamefont {A.~K.}\ \bibnamefont {Majee}},
  \bibinfo {author} {\bibfnamefont {A.}~\bibnamefont {Prohl}}, \ and\ \bibinfo
  {author} {\bibfnamefont {G.}~\bibnamefont {Vallet}},\ }\bibfield  {title}
  {\enquote {\bibinfo {title} {On stochastic optimal control in
  ferromagnetism},}\ }\href {\doibase 10.1007/s00205-019-01381-w} {\bibfield
  {journal} {\bibinfo  {journal} {Arch. Ration. Mech. Anal.}\ }\textbf
  {\bibinfo {volume} {233}},\ \bibinfo {pages} {1383--1440} (\bibinfo {year}
  {2019})}\BibitemShut {NoStop}%
\bibitem [{\citenamefont {Dunst}\ and\ \citenamefont
  {Prohl}(2016)}]{dunst2016}%
  \BibitemOpen
  \bibfield  {author} {\bibinfo {author} {\bibfnamefont {T.}~\bibnamefont
  {Dunst}}\ and\ \bibinfo {author} {\bibfnamefont {A.}~\bibnamefont {Prohl}},\
  }\bibfield  {title} {\enquote {\bibinfo {title} {The forward-backward
  stochastic heat equation: {N}umerical analysis and simulation},}\ }\href@noop
  {} {\bibfield  {journal} {\bibinfo  {journal} {SIAM J. Sci. Comput.}\
  }\textbf {\bibinfo {volume} {38}},\ \bibinfo {pages} {A2725--A2755} (\bibinfo
  {year} {2016})}\BibitemShut {NoStop}%
\bibitem [{\citenamefont {E}, \citenamefont {Han},\ and\ \citenamefont
  {Jentzen}(2017)}]{e2017}%
  \BibitemOpen
  \bibfield  {author} {\bibinfo {author} {\bibfnamefont {W.}~\bibnamefont {E}},
  \bibinfo {author} {\bibfnamefont {J.}~\bibnamefont {Han}}, \ and\ \bibinfo
  {author} {\bibfnamefont {A.}~\bibnamefont {Jentzen}},\ }\bibfield  {title}
  {\enquote {\bibinfo {title} {Deep learning-based numerical methods for
  high-dimensional parabolic partial differential equations and backward
  stochastic differential equations},}\ }\href@noop {} {\bibfield  {journal}
  {\bibinfo  {journal} {Commun. Math. Stat.}\ }\textbf {\bibinfo {volume}
  {5}},\ \bibinfo {pages} {349--380} (\bibinfo {year} {2017})}\BibitemShut
  {NoStop}%
\bibitem [{\citenamefont {Gorodetsky}, \citenamefont {Karaman},\ and\
  \citenamefont {Marzouk}(2018)}]{gorodetsky2018}%
  \BibitemOpen
  \bibfield  {author} {\bibinfo {author} {\bibfnamefont {A.}~\bibnamefont
  {Gorodetsky}}, \bibinfo {author} {\bibfnamefont {S.}~\bibnamefont {Karaman}},
  \ and\ \bibinfo {author} {\bibfnamefont {Y.}~\bibnamefont {Marzouk}},\
  }\bibfield  {title} {\enquote {\bibinfo {title} {High-dimensional stochastic
  optimal control using continuous tensor decompositions},}\ }\href@noop {}
  {\bibfield  {journal} {\bibinfo  {journal} {Int. J. Robot. Res.}\ }\textbf
  {\bibinfo {volume} {37}},\ \bibinfo {pages} {340--377} (\bibinfo {year}
  {2018})}\BibitemShut {NoStop}%
\bibitem [{\citenamefont {Kalise}\ and\ \citenamefont
  {Kunisch}(2018)}]{kalise2018}%
  \BibitemOpen
  \bibfield  {author} {\bibinfo {author} {\bibfnamefont {D.}~\bibnamefont
  {Kalise}}\ and\ \bibinfo {author} {\bibfnamefont {K.}~\bibnamefont
  {Kunisch}},\ }\bibfield  {title} {\enquote {\bibinfo {title} {Polynomial
  approximation of high-dimensional {H}amilton-{J}acobi-{B}ellman equations and
  applications to feedback control of semilinear parabolic {PDE}s},}\
  }\href@noop {} {\bibfield  {journal} {\bibinfo  {journal} {SIAM J. Sci.
  Comput.}\ }\textbf {\bibinfo {volume} {40}},\ \bibinfo {pages} {A629--A652}
  (\bibinfo {year} {2018})}\BibitemShut {NoStop}%
\bibitem [{\citenamefont {N\"{u}sken}\ and\ \citenamefont
  {Richter}(2021)}]{nuesken2021}%
  \BibitemOpen
  \bibfield  {author} {\bibinfo {author} {\bibfnamefont {N.}~\bibnamefont
  {N\"{u}sken}}\ and\ \bibinfo {author} {\bibfnamefont {L.}~\bibnamefont
  {Richter}},\ }\bibfield  {title} {\enquote {\bibinfo {title} {Solving
  high-dimensional {H}amilton-{J}acobi-{B}ellman {PDE}s using neural networks:
  perspectives from the theory of controlled diffusions and measures on path
  space},}\ }\href@noop {} {\bibfield  {journal} {\bibinfo  {journal} {Partial
  Differ. Equ. Appl.}\ }\textbf {\bibinfo {volume} {2}} (\bibinfo {year}
  {2021})}\BibitemShut {NoStop}%
\bibitem [{\citenamefont {Oster}, \citenamefont {Sallandt},\ and\ \citenamefont
  {Schneider}(2022)}]{oster2022}%
  \BibitemOpen
  \bibfield  {author} {\bibinfo {author} {\bibfnamefont {M.}~\bibnamefont
  {Oster}}, \bibinfo {author} {\bibfnamefont {L.}~\bibnamefont {Sallandt}}, \
  and\ \bibinfo {author} {\bibfnamefont {R.}~\bibnamefont {Schneider}},\
  }\bibfield  {title} {\enquote {\bibinfo {title} {Approximating optimal
  feedback controllers of finite horizon control problems using hierarchical
  tensor formats},}\ }\href@noop {} {\bibfield  {journal} {\bibinfo  {journal}
  {SIAM J. Sci. Comput.}\ }\textbf {\bibinfo {volume} {44}} (\bibinfo {year}
  {2022})}\BibitemShut {NoStop}%
\bibitem [{\citenamefont {Richter}, \citenamefont {Sallandt},\ and\
  \citenamefont {N\"usken}(2021)}]{richter2021}%
  \BibitemOpen
  \bibfield  {author} {\bibinfo {author} {\bibfnamefont {L.}~\bibnamefont
  {Richter}}, \bibinfo {author} {\bibfnamefont {L.}~\bibnamefont {Sallandt}}, \
  and\ \bibinfo {author} {\bibfnamefont {N.}~\bibnamefont {N\"usken}},\
  }\href@noop {} {\enquote {\bibinfo {title} {Solving high-dimensional
  parabolic pdes using the tensor train format},}\ }\bibinfo {howpublished}
  {preprint, https://arxiv.org/abs/2102.11830} (\bibinfo {year}
  {2021})\BibitemShut {NoStop}%
\bibitem [{\citenamefont {Sirignano}\ and\ \citenamefont
  {Spiliopoulos}(2018)}]{sirignano2018}%
  \BibitemOpen
  \bibfield  {author} {\bibinfo {author} {\bibfnamefont {J.}~\bibnamefont
  {Sirignano}}\ and\ \bibinfo {author} {\bibfnamefont {K.}~\bibnamefont
  {Spiliopoulos}},\ }\bibfield  {title} {\enquote {\bibinfo {title} {D{G}{M}: A
  deep learning algorithm for solving partial differential equations},}\
  }\href@noop {} {\bibfield  {journal} {\bibinfo  {journal} {J. Comput. Phys.}\
  }\textbf {\bibinfo {volume} {375}},\ \bibinfo {pages} {1339--1364} (\bibinfo
  {year} {2018})}\BibitemShut {NoStop}%
\bibitem [{\citenamefont {Stannat}\ and\ \citenamefont
  {Wessels}(2021{\natexlab{b}})}]{stannat20212}%
  \BibitemOpen
  \bibfield  {author} {\bibinfo {author} {\bibfnamefont {W.}~\bibnamefont
  {Stannat}}\ and\ \bibinfo {author} {\bibfnamefont {L.}~\bibnamefont
  {Wessels}},\ }\bibfield  {title} {\enquote {\bibinfo {title} {Deterministic
  control of stochastic reaction-diffusion equations},}\ }\href@noop {}
  {\bibfield  {journal} {\bibinfo  {journal} {Evol. Equ. Control Theory}\
  }\textbf {\bibinfo {volume} {10}},\ \bibinfo {pages} {701--722} (\bibinfo
  {year} {2021}{\natexlab{b}})}\BibitemShut {NoStop}%
\bibitem [{\citenamefont {Stannat}\ and\ \citenamefont
  {Vogler}(2023)}]{stannat2023}%
  \BibitemOpen
  \bibfield  {author} {\bibinfo {author} {\bibfnamefont {W.}~\bibnamefont
  {Stannat}}\ and\ \bibinfo {author} {\bibfnamefont {A.}~\bibnamefont
  {Vogler}},\ }\href@noop {} {\enquote {\bibinfo {title} {Approximation of
  optimal feedback controls for stochastic reaction-diffusion equations},}\
  }\bibinfo {howpublished} {preprint, arXiv, submit/4710581} (\bibinfo {year}
  {2023+})\BibitemShut {NoStop}%
\bibitem [{\citenamefont {Tudor}(1990)}]{tudor1990}%
  \BibitemOpen
  \bibfield  {author} {\bibinfo {author} {\bibfnamefont {C.}~\bibnamefont
  {Tudor}},\ }\bibfield  {title} {\enquote {\bibinfo {title} {Quadratic control
  for stochastic systems defined by evolution operators and square integrable
  martingales},}\ }\href@noop {} {\bibfield  {journal} {\bibinfo  {journal}
  {Math. Nachr.}\ }\textbf {\bibinfo {volume} {147}},\ \bibinfo {pages}
  {205--218} (\bibinfo {year} {1990})}\BibitemShut {NoStop}%
\bibitem [{\citenamefont {Vogler}(2023)}]{vogler2023}%
  \BibitemOpen
  \bibfield  {author} {\bibinfo {author} {\bibfnamefont {A.}~\bibnamefont
  {Vogler}},\ }\href@noop {} {\enquote {\bibinfo {title} {Sfb910 feedback},}\
  }\bibinfo {howpublished} {Github.
  \href{https://github.com/AVoglerTu/SFB910Feedback}{https://github.com/AVoglerTu/SFB910Feedback}}
  (\bibinfo {year} {2023})\BibitemShut {NoStop}%
\bibitem [{\citenamefont {Pinkus}(1999)}]{Pinkus1999}%
  \BibitemOpen
  \bibfield  {author} {\bibinfo {author} {\bibfnamefont {A.}~\bibnamefont
  {Pinkus}},\ }\bibfield  {title} {\enquote {\bibinfo {title} {Approximation
  theory of the mlp model in neural networks},}\ }\href@noop {} {\bibfield
  {journal} {\bibinfo  {journal} {Acta Numerica}\ }\textbf {\bibinfo {volume}
  {8}},\ \bibinfo {pages} {143--195} (\bibinfo {year} {1999})}\BibitemShut
  {NoStop}%
\bibitem [{\citenamefont {Carmona}\ and\ \citenamefont
  {Lauri\`ere}(2022)}]{carmona2022}%
  \BibitemOpen
  \bibfield  {author} {\bibinfo {author} {\bibfnamefont {R.}~\bibnamefont
  {Carmona}}\ and\ \bibinfo {author} {\bibfnamefont {M.}~\bibnamefont
  {Lauri\`ere}},\ }\bibfield  {title} {\enquote {\bibinfo {title} {Convergence
  analysis of machine learning algorithms for the numerical solution of mean
  field control and games: {II}---the finite horizon case},}\ }\href@noop {}
  {\bibfield  {journal} {\bibinfo  {journal} {Ann. Appl. Probab.}\ }\textbf
  {\bibinfo {volume} {32}},\ \bibinfo {pages} {4065--4105} (\bibinfo {year}
  {2022})}\BibitemShut {NoStop}%
\end{thebibliography}%

\end{document}